\documentclass[reqno]{amsart}

\usepackage{amsmath,amsthm,amssymb,color,graphicx}
\usepackage{amscd,verbatim}
\usepackage[all]{xy}

\theoremstyle{plain}
\newtheorem{theorem}{Theorem}
\newtheorem{proposition}[theorem]{Proposition}
\newtheorem{lemma}[theorem]{Lemma}
\newtheorem{corollary}[theorem]{Corollary}

\theoremstyle{definition}
\newtheorem{definition}[theorem]{Definition}

\theoremstyle{remark}
\newtheorem*{remark}{Remark}

\newtheorem*{example}{Example}

\numberwithin{equation}{section}
\numberwithin{theorem}{section}
\numberwithin{figure}{section}
\numberwithin{table}{section}

\newcommand{\cA}{{\mathcal A}}
\newcommand{\cB}{{\mathcal B}}

\newcommand{\cD}{{\mathcal D}}

\newcommand{\cF}{{\mathcal F}}

\newcommand{\cK}{{\mathcal K}}
\newcommand{\cL}{{\mathcal L}}
\newcommand{\cM}{{\mathcal M}}

\newcommand{\cR}{{\mathcal R}}
\newcommand{\cS}{{\mathcal S}}

\newcommand{\cX}{{\mathcal X}}
\newcommand{\cY}{{\mathcal Y}}
\newcommand{\cZ}{{\mathcal Z}}

\newcommand{\alg}{{\cA}}
\newcommand{\cb}{{\cB}}
\newcommand{\cc}{{\cC}}
\newcommand{\hilb}{{\cH}}

\newcommand{\cpt}{{\cK}}

\newcommand{\CC}{{\mathbb C}}

\newcommand{\OO}{{\mathbb O}}

\newcommand{\RR}{{\mathbb R}}
\newcommand{\TT}{{\mathbb T}}

\newcommand{\ZZ}{{\mathbb Z}}

\newcommand{\id}{{\mathbf 1}}

\newcommand{\txtfrac}[2]{{\textstyle{{\frac{#1}{#2}}}}}
\newcommand{\wM}{{\widetilde M}}
\newcommand{\wH}{{\widetilde H}}
\newcommand{\wP}{{\widetilde P}}
\newcommand{\wGamma}{{\widetilde \Gamma}}
\newcommand{\PU}{{\sf PU}}
\newcommand{\PUH}{{\sf PU}({\mathcal H})}
\newcommand{\U}{{\sf U}}

\newcommand{\what}{\widehat}

\newcommand{\integer}{\ZZ}

\newcommand{\sfG}{{\mathsf G}}
\newcommand{\sfH}{{\mathsf H}}
\newcommand{\sfK}{{\mathsf K}}
\newcommand{\sfN}{{\mathsf N}}
\newcommand{\sfT}{{\mathsf T}}

\def\aut{{\rm Aut}}
\def\id{{\rm id}}
\def\alg{{\mathcal A}}
\def\cat{{\mathcal C}}
\def\hilb{{\mathcal H}}
\def\cm{{\mathcal M}}
\def\cpt{{\mathcal K}}
\def\bdd{{\mathcal B}}
\def\cc{{\mathcal C}}
\def\ind{{\rm ind}}
\def\uind{\hbox {$u$-ind}}
\def\ad{{\rm ad}}
\def\idob{{\bf 1}}
\def\integer{{\bf Z}}

\def\wt#1{\widetilde{#1}}
\def\wh#1{\widehat{#1}}
\def\dual#1{\wh{#1}}
\def\ip#1#2{\langle{#1},{#2}\rangle}
\def\dip#1#2{\langle\kern-2pt\ip{#1}{#2}\kern-2pt\rangle}
\def\lip#1#2{\ip{#1}{#2}_L}
\def\rip#1#2{\ip{#1}{#2}_R}
\def\conj#1{\overline{#1}}

\def\sect#1{\advance\count30 by 1\count31=0 
\bigskip\noindent{\bf \the\count30. #1}

\bigskip} 
\def\dfn{\advance\count32 by 1
\bigskip\noindent{\bf Definition \the\count30.\the\count32.  }}
\def\thm#1{\advance\count31 by 1
\bigskip\noindent{\bf #1 \the\count30.\the\count31.  }}
\def\athm#1{\advance\count31 by 1
\bigskip\noindent{\bf #1 A.\the\count31.  }}

\begin{document}

\title[C$^*$-algebras in tensor categories]{C$^*$-algebras in tensor categories}

\author[P Bouwknegt]{Peter Bouwknegt}

\address[Peter Bouwknegt]{
Department of Theoretical Physics,
Research School of Physics and Engineering, and 
Department of Mathematics, Mathematical Sciences Institute, 
The Australian National University, 
Canberra, ACT 0200, Australia}
\email{peter.bouwknegt@anu.edu.au}

\author[KC Hannabuss]{Keith C. Hannabuss}

\address[Keith Hannabuss]{
Mathematical Institute, 24-29 St. Giles', Oxford, OX1 3LB, and 
Balliol College, Oxford, OX1 3BJ,
England}
\email{kch@balliol.oxford.ac.uk}

\author[V Mathai]{Varghese Mathai}

\address[Varghese Mathai]{
Department of Pure Mathematics,
University of Adelaide, 
Adelaide, SA 5005, 
Australia}
\email{mathai.varghese@adelaide.edu.au}

\begin{abstract}
We define and systematically study nonassociative C$^*$-algebras as
C$^*$-algebras internal to a topological tensor category. We also offer a 
concrete approach to these C$^*$-algebras, as $\sfG$-invariant, 
norm closed $*$-subalgebras of
bounded operators on a ${\sfG}$-Hilbert space, with deformed composition product. 
Our central results are those of stabilization and Takai duality for (twisted) crossed
products in this context. 
\end{abstract}



\maketitle

\section{Introduction} \label{secA}

In \cite{BHM3} we gave an account of some nonassociative algebras and their 
applications to T-duality, with a brief mention of the role of categories at 
the end.
In this paper we will develop the theory more systematically from the category 
theoretic perspective.
In particular we shall not need to assume that the groups are abelian, and 
will mostly work with general three-cocycles rather than antisymmetric 
tricharacters. We believe, however, that the results in
this paper are of interest independent of our original motivation.
Since writing \cite{BHM3} we have become aware of more of the large literature in 
this subject, for example, the work of Fr\"ohlich, Fuchs, Runkel, Schweigert 
in conformal field theory \cite{FFRS1, FFRS2, FRS1, FRS2, FRS3, FRS4}, 
and of Beggs, Majid and collaborators  \cite{AM1, AM2, BM1, BM2, BM3, BM4}.
Most of that work is algebraic in spirit, working with finite groups or 
finite dimensional Hopf algebras, whereas we are primarily interested in 
locally compact groups and C$^*$-algebras, which necessitate the development 
of a rather different set of techniques.
The work of Nesterov and collaborators \cite{N1, N2, N3, N4, N5} does use vector groups, but has 
rather different aims and methods.
(In addition, there is a well-established theory of C$^*$- and W$^*$-categories, 
\cite{GLR, M}, in which the algebras are generalised to morphisms in a suitable 
category, which means that they are automatically associative.
One could generalise to weak higher categories, but the examples discussed in \cite{BHM3} 
suggest that
one starts by looking at categories in which the algebras and their modules are objects, 
which is also more directly parallel to the algebraic cases already mentioned.

In Section \ref{secB} 
we introduce the tensor categories that we use and give some elementary examples, 
based on our earlier work in \cite{BHM3}.
Within the category it is possible to define algebras and modules.
The next three sections show how to obtain nonassociative algebras of twisted 
compact and bounded operators on a Hilbert space, and introduce Morita 
equivalence.
We then link this to exterior equivalence in Section \ref{secF}, which establishes that
exterior equivalent twisted actions give rise to isomorphic twisted crossed product
C$^*$-algebras. 

Section \ref{secG} is devoted to an extension of the nonassociative Takai duality 
proved in \cite{BHM3}.
This is especially useful, because it provides a method of stabilising 
algebras.
For example, an associative algebra with a very twisted group action has a
nonassociative dual and double dual which admit ordinary untwisted group 
actions.
The double dual is Morita equivalent to the original, and so one could replace 
the original associative algebra and twisted action by an equivalent 
nonassociative algebra and ordinary action.

The main result in Section \ref{secH} establishes the fact that twisted crossed products 
can be obtained by repeated ordinary crossed products, but with a possible
modified automorphism action of the final subgroup, a result that goes a long 
way towards proving an analog of the Connes-Thom isomorphism theorem 
\cite{Co81, Co94} in our context, as briefly discussed in the final section.

A theorem of MacLane \cite{M} asserts that every monoidal category can be made 
strict, that is, associative, but in general the functor which does this is 
quite complicated. 
However, in our category things are much simpler, and in Appendix \ref{secappA}, it is 
shown that whenever a nonassociative algebra acts on a module, its 
multiplication can be modified to an associative multiplication.
Examples of the strictification process are discussed there,
in particular to the algebras of twisted compact and bounded operators which 
are defined to act on a module, but also have more serious implications for 
physics, where the algebras are generally represented by actions on modules.
(On the other hand this does not mean that we can simply dismiss the 
nonassociativity, because there are known nonassociative algebras such as the 
octonions, cf. \cite{Baez}, which fit into our framework.) Appendix \ref{secappB} gives a concrete approach
to our nonassociative C$^*$-algebras, as $\sfG$-invariant, 
norm closed $*$-subalgebras of bounded operators on a ${\sfG}$-Hilbert space, 
with composition product deformed by a 3-cocycle on $\sfG$. In Appendix \ref{secappC}, 
we revisit the construction of our nonassociative torus, via a geometric construction
that realizes it as a nonassociative deformation of the C$^*$-algebra of  
continuous functions on the torus. For completeness, we have summarized in Appendix \ref{secappD}, 
the original motivation for this 
work, namely, T-duality in string theory.

\medskip

{\bf Acknowledgments.} KCH would like to thank E.~Beggs and S.~Majid for discussions 
at the ``Noncommutative Geometry'' program held at the Newton Institute (Cambridge)
in 2006. VM would like to thank A.~Connes for feedback and suggestions at the 
``Noncommutative Geometry'' workshop held at Banff, Canada in 2006, and PB
would like to thank the participants at the ``Morgan-Phoa mathematics workshop'' 
at the ANU for useful feedback.
Part of this work was done while we were visiting the Erwin Schr\"odinger Institute
in Vienna for the program on ``Gerbes, Groupoids and Quantum Field Theory''.
KCH would also like to thank the University of Adelaide for hospitality 
during the intermediate stages of the project.
Both PB and VM acknowledge financial support from the Australian Research Council.

\section{Tensor categories and their algebras} \label{secB}

As in \cite{K} we use tensor category to  mean just a monoidal category, without 
any of the other structures often assumed elsewhere.
That is, a tensor category is a category in which associated to each 
pair of objects $A$ and $B$ there exists 
a product object $A\otimes B$, and there is an identity object $\idob$, 
such that $\idob\otimes A \cong A \cong A\otimes \idob$,
together with associator isomorphisms
\begin{equation}
\Phi = \Phi_{A,B,C}: A\otimes (B\otimes C) \to (A\otimes B)\otimes C
\end{equation}
for any three objects $A$, $B$ and $C$, satisfying the consistency pentagonal identity on quadruple products:

\begin{equation*}
\xymatrix@C-50pt{
&& A\otimes (B\otimes(C\otimes D)) \ar[ddll]_{\idob_A\otimes \Phi_{B,C,D}}  
\ar[ddrr]^{\Phi_{A,B,C\otimes D}} && \\ &&&& \\
A\otimes ((B\otimes C))\otimes D) \ar[dddr]^{\Phi_{A,B\otimes C,D}} && &&
(A\otimes B)\otimes (C\otimes D) \ar[dddl]_{\Phi_{A\otimes B,C,D}}  \\ &&&& \\ &&&& \\
& (A\otimes (B\otimes C))\otimes D \ar[rr]_{\Phi_{A,B,C} \otimes \idob_D} && 
((A\otimes B)\otimes C))\otimes D & }
\end{equation*}
with each arrow the appropriate map $\Phi$, and the triangle relation:
\begin{equation*}
\xymatrix{
A\otimes(\idob\otimes B) \ar[ddr]_{\cong} \ar[rr]_{\Phi_{A,\idob,B}} && (A\otimes \idob)\otimes B 
\ar[ddl]^{\cong} \\ &&\\
& A\otimes B   }  
\end{equation*}
MacLane's coherence theorem ensures that these conditions are sufficient to 
guarantee consistency of all other rebracketings, \cite{McL}.
(The theorem proceeds by showing that one can always take the category to be a strictly 
associative category.  We prove this explicitly for our examples in Appendix \ref{secappA})

Module categories provide two standard examples of tensor categories with the obvious 
identification map 
$\Phi = \id$.
One is the category of ${\cR}$-modules and  ${\cR}$-morphisms, 
for ${\cR}$ a commutative algebra over $\CC$.
The appropriate tensor product of objects $A$ and $B$ is $A\otimes_{\cR}B$, and the identity 
object is ${\cR}$ itself. 
Continuous trace algebras with spectrum $S$ are $C_0(S)$-modules, and so can be studied 
within this 
category with ${\cR} = C_0(S)$, though some care is needed in defining the appropriate tensor 
products in 
the case of topological algebras, but this can be done explicitly in this case, 
cf. \cite[Sect 6.1]{RW}.

A subtly different example is provided by the algebra of functions $H = C_0(\sfG)$ on 
a separable locally compact group $\sfG$.
As well as being an algebra under pointwise multiplication it also has a comultiplication 
$\Delta: H \to H\otimes H$ taking a function $f\in C_0(\sfG)$ to $(\Delta f)(x,y) = f(xy)$, and a counit 
$\epsilon: H \to \CC$, which evaluates $f$ at the group identity.
This enables us to equip the category of $H$-modules with a tensor product $A\otimes B$ over 
$\CC$, on which $f$ acts as $\Delta(f)$, 
and the identity object being $\CC$ with the  trivial $H$-action of multiplication by $\epsilon(f)$.
Writing the comultiplication in abbreviated Sweedler notation $\Delta f = f_{(1)}\otimes f_{(2)}$,
the action on $A\otimes B$ is
$$f[a\star b] = f_{(1)}[a]\star f_{(2)}[b].$$
Since everything has been defined in terms of the comultiplication and counit 
of $H$, this clearly generalises to bialgebras, and even to quasi-bialgebras.
If $\sfG$ is an abelian group,
with Pontryagin dual group $\dual{\sfG} = \text{Hom}(\sfG, \mathsf{U}(1))$, 
then we can work with  $C^*(\dual{\sfG})$ 
instead of $C_0(\sfG)$, and the tensor product 
action of $\xi\in \dual{\sfG}$ is just $\xi\otimes \xi$. 
If the group $\dual{\sfG}$ acts on $S$, the two examples can be 
combined in the category of modules for the crossed product, or transformation 
groupoid, algebra $C_0(S)\rtimes \dual{\sfG}$, equipped with the tensor product over 
$C_0(S)$, and identity object $C_0(S)$.

Both examples use modules for a $*$-algebra, with $f^*(s) = \conj{f(s)}$ in $C_0(S)$ and 
$f^*(x) = \conj{f(x^{-1})}$ in $C_0({\sfG})$, and the category contains conjugate objects $A^*$, 
having the same underlying set, but with the algebra action changed to that of $f^*$ and 
conjugated scalar multiplication. 
For $a$ an element of some object $A$ and $a^*$ the same element considered as an element of 
$A^*$ we then have $f^*[a^*] = f[a]^*$.
Working with $C^*(\dual{{\sfG}})$ instead of $C_0({\sfG})$ one has $\conj{\xi(x^{-1})} = \xi(x)$ so that 
$\xi^* = \xi$.
This too can be generalised to quasi-Hopf algebras \cite[Section XV.5]{K} which  
provide such structure as do the coinvolutions in Kac C$^*$-algebras.
Conjugation $A \mapsto A^*$ preserves direct sums and gives a covariant functor. 
Another crucial property follows by noting that in $C_0({\sfG})$ one has
$$\Delta(f^*)(x,y) = f^*(xy) = \conj{f(y^{-1}x^{-1})} =
\conj{\Delta(f))(y^{-1},x^{-1})} =  (f_{(2)}^*\otimes f_{(1)}^*)(x,y),$$
so that conjugation reverses the order of factors in the tensor product.
There is a natural isomorphism  between $(A\otimes B)^*$ and $B^*\otimes A^*$.
Since preparing the first draft of this paper the preprint \cite{BM3} has 
appeared and gives a systematic account of bar categories, of which these form 
one example. We refer the reader there for more detail.

Tensor categories have enough structure to define algebras and modules.

\begin{definition}\label{defAa}
An object $\alg$ is an algebra (or monoid) in a tensor category if there is a 
morphism $\star: \alg \otimes \alg \to \alg$ that is associative in the category, that is,  
$\star(\star\otimes\id)\Phi = \star(\id \otimes\star)$ as maps $\alg\otimes (\alg\otimes \alg) \to \alg$.
An algebra $\alg$ in the category is called a $*$-algebra if 
the category has the above conjugation of objects and 
$\alg^* = \alg$.
A left (respectively, right) module $\cM$ for $\alg$ is an object such that 
there is a morphism, which we also denote by $\star$, sending 
$\alg\otimes \cM$ to $\cM$ (respectively, $\cM\otimes \alg$ to $\cM$), and satisfying 
the usual composition law in the category, that is, for left modules 
$\star(\star\otimes\id)\Phi = \star(\id\times\star)$ as maps $\alg\otimes (\alg\otimes \cM) \to \cM$.
For brevity, the term module will mean a left module, unless otherwise 
specified.
An $\alg$-$\bdd$-bimodule $\cX$ for two algebras $\alg$ and $\bdd$ in the 
category is a left module for $\alg$ and a right module for $\bdd$, with 
commuting actions (allowing for rebracketing given by an associator map).
\end{definition}

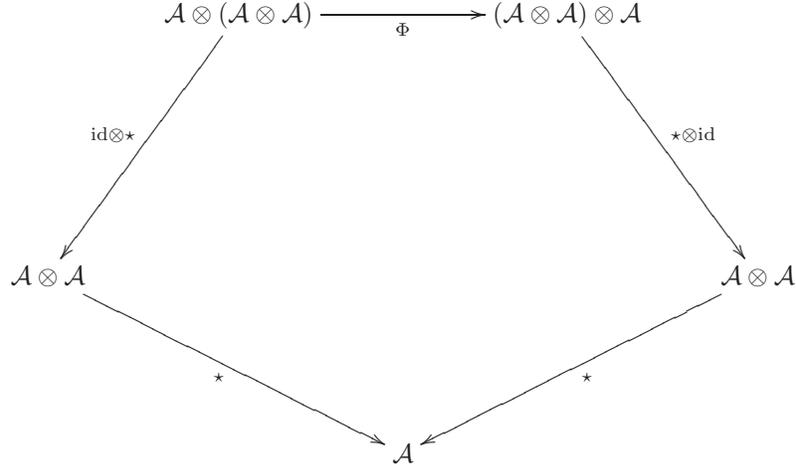
\begin{equation*}
\xymatrix{
&\cA \otimes ( \cA \otimes \cA ) \ar[dddl]_{\id\otimes \star}  \ar[rr]_\Phi && 
(\cA \otimes  \cA) \otimes \cA \ar[dddr]^{\star\otimes\id} & \\
&&&& \\ &&&& \\
\cA \otimes \cA \ar[ddrr]_\star &&&& \cA \otimes \cA \ar[ddll]^\star \\
&&&& \\
&& \cA }
\end{equation*}

As an example we note that continuous trace algebras with spectrum $S$ are algebras in the category  of 
$C_0(S)$-modules introduced above.
In fact $C_0(S)$ can be identified with a subalgebra of the centre $ZM(\alg)$ of the multiplier algebra of a 
continuous trace algebra $\alg$, and so it also acts on all $\alg$-modules, so that $\alg$-modules in the 
ordinary sense are also $\alg$-modules in the category.
Similarly, continuous trace algebras with spectrum $S$ on which $\dual{{\sfG}}$ acts as automorphisms are 
algebras in the category of $C_0(S)\rtimes \dual{{\sfG}}$-modules.
Although we have so far taken $\Phi$ to be the usual identification map, there are other possibilities.
Let ${\sfG}$ be a separable locally compact group with dual $\wh{{\sfG}}$, and 
let $\phi\in C({\sfG}\times {\sfG}\times {\sfG})$ be normalised to take the value 1 whenever 
any of its arguments is the identity $1\in {\sfG}$, and satisfy the pentagonal 
cocycle identity
\begin{equation}
\phi(x,y,z)\phi(x,yz,w)\phi(y,z,w) = \phi(xy,z,w)\phi(x,y,zw) \,.
\end{equation}
Since $C({\sfG})$ is the multiplier algebra of $C_0({\sfG})$ it also acts on $C_0({\sfG})$-modules.
Unlike the algebraic case there are various module tensor products, and it 
is assumed that we have chosen one for which the action of 
$C({\sfG})\otimes C({\sfG})\otimes C({\sfG})$ is defined.

\begin{definition}\label{defAb}
The category $\cat_{\sfG}(\phi)$ has for objects normed $C_0({\sfG})$-modules, or 
equivalently normed $\wh{{\sfG}}$-modules, and its morphisms are continuous linear maps commuting 
with the action. 
The tensor  structure comes from taking the tensor product of modules 
with the tensor product action of the coproduct $(\Delta f)(x,y) = f(xy)$ for 
$f\in C({\sfG})$, (or diagonal tensor product action $\xi\otimes\xi$ 
of  $\xi\in\wh{{\sfG}}$). 
For any three objects $\alg$, $\bdd$ and $\cc$ the associator map, 
$\Phi: \alg\otimes (\bdd\otimes \cc) \to (\alg\otimes \bdd)\otimes \cc$ is 
given by the action of 
$\phi\in C({\sfG}\times {\sfG}\times {\sfG}) \cong C({\sfG})\otimes C({\sfG})\otimes C({\sfG})$. 
The identity object is the trivial one-dimensional module $\CC$, on which 
$\wh{{\sfG}}$ acts trivially, or, equivalently, $f\in C({\sfG})$ multiplies by $f(1)$, 
where 1 is the identity element in ${\sfG}$. 
\end{definition}

We could introduce a similar structure for $C_0(S)\rtimes \dual{{\sfG}}$ modules.

As already mentioned, the algebras (or monoids) in $\cat_{\sfG}(\phi)$ are objects $\alg$ 
for which there is a product morphism $\alg\otimes \alg \to \alg$, which we shall write 
as $a\otimes b\mapsto a\star b$. 
An algebra $\alg$ can have a module $\cm$, when there is a morphism 
$\alg\otimes \cm \to \cm$, which we also write $a\otimes m \mapsto a\star m$, relying on 
the context to distinguish multiplications from actions.
The category structure forces interesting compatibility conditions on algebras and modules.

\begin{proposition} \label{thmAc}
Let $\alg$ be an algebra in the category $\cat_{\sfG}(\phi)$.
Then the group $\dual{{\sfG}}$ acts on an algebra $\alg$ by automorphisms, and the action of 
$\dual{{\sfG}}$ on an $\alg$-module $M$ gives  a covariant representation of $\alg$ and $\dual{{\sfG}}$
or, equivalently, 
a representation of the crossed product $\alg\rtimes \dual{{\sfG}}$. 
\end{proposition}

\begin{proof}
For $a, b\in \cA$, and $\xi\in \dual{{\sfG}}$ we have $\xi[a]\star\xi[b] = \xi[a\star b]$, 
showing that the action of $\dual{{\sfG}}$ is by endomorphisms, and since $\dual{{\sfG}}$ 
is a group, these are invertible, and so automorphisms. 

The morphism property gives $\xi[a]\star\xi[m] = \xi[a\star m]$, showing that the actions 
of $\alg$ and $\dual{{\sfG}}$ combine into a covariant representation of 
$(\alg,\dual{{\sfG}})$.  Standard theory then tells us that this is equivalent to having an action of 
the crossed product $\alg\rtimes \dual{{\sfG}}$. 
\end{proof}

We need to take care concerning the ordering of products where $\Phi$ sets up the appropriate 
associativity conditions.
When $\alg$ is an algebra in the category, and $\cM$ is an $\alg$-module,
we shall simplify the notation by writing 
$\Phi(a\star (b\star m))$ for the composition of the maps 
$\Phi: \alg\otimes(\alg\otimes \cM) \to (\alg\otimes\alg)\otimes \cM$ and the multiplications and action
\begin{equation*}
(\alg\otimes\alg)\otimes \cM \to  \alg\otimes \cM \to \cM\,,
\end{equation*}
applied to $a\otimes(b\otimes c)$, and similarly for triple products of algebra elements.
We shall similarly abbreviate the notation for other operations involving tensor products.

The algebra of twisted compact operators $\cpt_\phi(L^2(G))$, introduced in 
\cite{BHM3}, provides an example of an algebra in the category $\cat_{\sfG}(\phi)$.
There we assumed that $\phi$ was an antisymmetric tricharacter, that is,
$\phi(x,y,z)$ is a character in each of its arguments for fixed values of the others, and is 
inverted by any transposition of its arguments.
This is sufficient to ensure that it satisfies the cocycle identity, but at the expense of slightly more 
complicated formulae the cocycle condition usually suffices, for example the 
product of two integral kernels is defined by
\begin{equation}
(k_1\star k_2)(x,z) 
= \int \phi(xy^{-1},yz^{-1},z)k_1(x,y)k_2(y,z)\,dydz \,.
\end{equation}
The action of $f\in C(\sfG)$ on $\cpt_\phi(L^2(\sfG))$ multiplies the kernel $K(x,y)$ by 
$f(xy^{-1})$, and one checks that this defines an automorphism
using the identity $(\Delta f)(xy^{-1},yz^{-1}) = f(xz^{-1})$.
Moreover, as may be readily checked using the pentagonal identity:
\begin{align}
&((k_1\star k_2)\star k_3))(x,w) \\
&= \int \phi(xy^{-1},yz^{-1},z)\phi(xz^{-1},zw^{-1},w)
k_1(x,y)k_2(y,z)k_3(z,w)\,dydz \nonumber \\
&= \int \phi(xy^{-1},yz^{-1},zw^{-1})
\phi(xy^{-1},yw^{-1},w)\phi(yz^{-1},zw^{-1},w) \nonumber \\
& \qquad \times k_1(x,y)k_2(y,z)k_3(z,w)\,dydz\nonumber \\
&=\Phi(k_1\star(k_2\star k_3))(x,w) \,.
\end{align}

In Appendix \ref{secappB} the twisted compact  operators $\cpt_\phi(L^2({\sfG}))$, introduced in 
\cite{BHM3}, and the twisted bounded operators $\bdd_\phi(L^2({\sfG}))$
are systematically defined and studied, providing examples of C$^*$-algebras in the category 
$\cat_{\sfG}(\phi)$. In fact, any norm closed, ${\sfG}$-invariant $\star$-subalgebra of $\bdd_\phi(L^2({\sfG}))$ 
gives  an example of a C$^*$-algebra in the category. The reworking of the twisted compact 
operators and  twisted bounded operators with a general cocycle $\phi$ 
shows that the arguments of \cite{BHM3} need not depend on $\phi$ being a tricharacter, 
although there may be simplifications when it is.

\begin{proposition} \label{thmAd}
Let $\lambda: A\otimes(B\otimes C) \to \CC$ be a morphism in $\cat_{\sfG}(\phi)$.
If $\phi$ is a tricharacter then 
$\lambda(\Phi(a\otimes (b\otimes c))) = \lambda(a\otimes (b\otimes c))$, 
for all $a\in A$, $b\in B$ and $c\in C$.
The same applies for a morphism $\lambda:A\otimes (B\otimes C) \to A\otimes \CC \to A$
which factors through $\CC$.
\end{proposition}

\begin{proof}
The fact that $\lambda$ is a morphism means that for any $f\in C_0({\sfG})$
\begin{equation*}
\epsilon(f)\lambda\Phi = \lambda(\id\otimes \Delta)\Delta(f))\Phi
=  \lambda(f_{(1)}\otimes (f_{(2)}\otimes f_{(3)}))\phi \,,
\end{equation*}
or 
\begin{equation*}
\lambda\epsilon(f)\Phi = \lambda(f_{(1)}\otimes (f_{(2)}\otimes f_{(3)}))\phi \,,
\end{equation*}
Thus the effect of $(f_{(1)}\otimes (f_{(2)}\otimes f_{(3)}))\phi(x,y,z) = f(x(yz))\phi(x,y,z)$ 
is just the same as $\epsilon(f)\phi(x,y,z) = f(1)\phi(x,y,z)$.
In other words, the effect of $\lambda\Phi$ is concentrated where $xyz = 1$.
However, when $\phi$ is an antisymmetric  tricharacter, it takes the value 1 when arguments are repeated, 
and so  
$$\phi(x,y,z) = \phi(x,y,x)\phi(x,y,y)\phi(x,y,z)= \phi(x,y,xyz) = 1.$$
In other words $\lambda\Phi = \lambda$. 
A similar argument applied to $B\otimes C$
covers the case when $\lambda$ factors through $\CC$.
 
We conclude by remarking that in the abbreviated notation introduced above,
one would write
$\lambda(\Phi(a\otimes (b\otimes c))) = \lambda(a\otimes (b\otimes c))$.
\end{proof}

More generally $C({\sfG})$ can be replaced by other algebras.
To give the desired structure we require at least a comultiplication 
$\Delta$ to define a tensor product action, a linear functional $\epsilon$ 
satisfying $\epsilon(h_1)h_2 = h = h_1\epsilon(h_2)$ defining the action on 
the identity object, and a three-cocycle $\Phi$ in the multiplier algebra
$\cm(\sfH\otimes \sfH\otimes \sfH)$, consistent with these.
The pentagonal cocyle condition is 
\begin{equation}
(\Delta\otimes \id\otimes \id)(\Phi )(\id\otimes\id\otimes\Delta)(\Phi)
= (\Phi\otimes \id)(\id\otimes\Delta\otimes\id)(\Phi)(\id\otimes \Phi) \,.
\end{equation}
Consistency of the associativity rebracketing with the action of $\sfH$ on 
tensor products of modules requires 
\begin{equation}
(\Delta\otimes\id)\Delta(h)= \ad\Phi(\id\otimes\Delta)\Delta(h) \,,
\end{equation}
whilst consistency with the action on the identity object means that 
$\epsilon$ contracted with the middle part of $\Phi$ gives the identity.
These are precisely the conditions satisfied by a quasi-bialgebra 
\cite[Section XV.1]{K}.

\begin{definition} \label{defAe}
The objects in the category $\cat^\sfH(\Phi)$ are $\sfH$-modules, and the morphisms are
linear $\sfH$-endomorphisms,  with the associator map given by the action of 
$\Phi\in \sfH\otimes \sfH\otimes \sfH$.
The action on tensor products is given by the comultiplication and the 
trivial object is $\CC$ with the action given by the counit.
\end{definition}

The quasi-bialgebra version of Proposition \ref{thmAc} asserts that 
$\sfH$ acts by automorphisms of an algebra $\alg$ in the category, that is, 
$h[a\star b] = h_{(1)}[a]\star h_{(2)}[b]$,
for all $a, b\in A$,  and the actions of $\sfH$ and $\alg$ on an $A$-module 
are covariant, so that $h[a\star m] = h_{(1)}[a]\star h_{(2)}[m]$,
for all $m\in \cm$, and one has a module for the crossed product $A\rtimes \sfH$.

We shall often use $\sfH$ to include the case of $C({\sfG})$, though in the latter 
case there are extra analytic conditions. 
(In principle one might use Kac algebras, \cite{ES, EV}, but that would require too 
big a digression.)

We shall outline our results for topological groups (more directly 
linked to our applications) and for quasi-Hopf algebras. 

We shall work in tensor categories of modules for an appropriate 
locally compact abelian group (more directly linked to our applications) or for 
a quasi-Hopf algebra (which makes the algebraic structure particularly 
transparent), though it is probably 
possible to extend much of this to tensor categories of modules for 
Kac algebras \cite{ES, EV}, which provide a natural framework for considerations 
of duality. 
Even locally compact groups present challenges beyond those present in 
the purely algebraic case of Hopf algebras, for example, the algebra $C_0({\sfG})$ 
of compactly supported functions on ${\sfG}$ has no unit, since the constant 
function 1 is not compactly supported, whilst $L^\infty({\sfG})$ has no counit since 
evaluation at the identity is not defined, and in neither case is there an 
antipode, though both have a coinvolution as Kac C$^*$-algebras.
One also requires modules which are not finitely generated.

\section{Hilbert modules in tensor categories} \label{secC}

The notion of the twisted kernels on $L^2({\sfG})$ can be extended to more general 
Hilbert spaces. 
There are several equivalent characterisations of the bounded operators in a normal Hilbert 
space, and for our purposes the most useful approach is rather indirect.

Let $\hilb$ be a $\wh{{\sfG}}$-module with an inner product $\ip{\cdot}{\cdot}$.
Since the inner product takes values in the trivial object $\CC$, 
consistency with the action of $\wh{{\sfG}}$ requires that
$\ip{\xi[\psi_1]}{\xi[\psi_2]} = \xi[\ip{\psi_1}{\psi_2}] = \ip{\psi_1}{\psi_2}$, 
for all $\xi\in \wh{{\sfG}}$ and 
$\psi_1$, $\psi_2 \in \hilb$, so that the $\wh{{\sfG}}$ action is unitary, or equivalently
$\ip{\psi_1}{f[\psi_2]} = \ip{f^*[\psi_1]}{\psi_2}$ for all $f\in C({\sfG})$.

\begin{definition}
An object $\hilb$ in $\cat_{\sfG}(\phi)$ with an inner product $\ip{\cdot }{\cdot }$ with respect to 
which the action of $\wh{{\sfG}}$ is unitary (or, equivalently, consistent with the 
$*$-structure of $C({\sfG})$) is called a pre-Hilbert space in $\cat_{\sfG}(\phi)$.
If it is complete in the norm topology it is called a Hilbert space in 
$\cat_{\sfG}(\phi)$.
\end{definition}

The unitarity of the action of $\wh{{\sfG}}$ means that the inner product on $\hilb$ 
defines a morphism  from the conjugate space $\hilb^*$ introduced earlier to the 
dual of $\hilb$ given by $\psi^* = \ip{\psi}{\cdot}$.
We can alternatively think of an inner product as a morphism $\hilb^*\otimes\hilb \to \CC$, 
written $\psi_1^*\otimes\psi_2 \mapsto \ip{\psi_1}{\psi_2}$, which satisfies the positivity condition.

\begin{proposition}
In order that the map taking $A$ to $A^*$ be consistent with the associator 
isomorphism it is necessary and sufficient that 
\begin{equation}
\Phi(A\otimes(B\otimes C))^* =  \Phi^{-1}((A\otimes (B\otimes C))^*) \,,
\end{equation}
for all objects $A$, $B$ and $C$, or, equivalently in $\cat_{\sfG}(\phi)$, that the $\phi$ be unitary. 
\end{proposition}

\begin{proof}
We have
\begin{align*}
\Phi(A\otimes(B\otimes C))^*
&= ((A\otimes B)\otimes C)^*\cr
&=  C^*\otimes (B^*\otimes A^*)\cr
&=  \Phi^{-1}((C^*\otimes B^*)\otimes A^*)\cr
&=  \Phi^{-1}((A\otimes (B\otimes C))^*),
\end{align*}
which amounts to saying that the function $\phi^* = \phi^{-1}$,
so that $\phi$ is unitary.
\end{proof}

In ordinary Hilbert spaces bounded operators can be characterised as those that are adjointable, and 
this definition is easy to generalise.

\begin{definition}
Modifying the usual definition we shall call a linear operator $A$ on $\hilb$ 
adjointable if for all $\psi_1$ and $\psi_2\in \hilb$ there is an operator 
$A^\star$ such that
$$\ip{A\star\psi_1}{\psi_2} = \Phi(\ip{\psi_1}{A^\star\star\psi_2}),$$ 
where the ordering on each side is that given by the order-reversing conjugation, that is, the two 
sides are images of  
$(\psi_1^*\otimes A^*)\otimes \psi_2$ and $\Phi(\psi_1^*\otimes(A^\star\otimes\psi_2))$ under the module 
action and inner product map.
\end{definition}

We note that this definition is taken to mean that $A$ and $A^\star$ are in a $*$-algebra with a 
module $\hilb$. 
The unitarity of the action means that the ordinary bounded operators are objects in the category 
with the action $\xi[A] = \xi\circ A\circ \xi^{-1}$.
The equality of the images of $(\psi_1^*\otimes A^*)\otimes \psi_2$ and  
$\Phi(\psi_1^*\otimes(A^\star\otimes\psi_2))$ allows us to identify the adjoint $A^\star$ with the conjugate $A^*$.
We note that the unique action of $\wh{{\sfG}}$ consistent with the 
covariance property gives $\xi: A \mapsto \xi\circ A\circ\xi^{-1}$.
When $\Phi$ is given by an antisymmetric tricharacter $\phi$ we may apply Proposition 
\ref{thmAd} to the 
morphism $\lambda(\psi_1\otimes(A\otimes\psi_2)) = \ip{\psi_1}{A\star\psi_2}$
to deduce that $\Phi$ acts trivially and the condition for the adjoint reduces to 
$\ip{A^\star\psi_1}{\psi_2} = \ip{\psi_1}{A\psi_2}$, as usual.
In this case the adjointable operators are therefore just the bounded operators on $\hilb$.
For general $\phi$ one will have a natural generalisation of the bounded operators, the subject of the next 
section.

\section{Twisted compact and twisted bounded operators} \label{secD}

Regarding $\hilb$ as a right $\CC$-module for the scalar 
multiplication action,  Rieffel's method allows us to define the dual 
inner product $\dip{\cdot}{\cdot}$ such that 
$$\dip{\psi_0}{\psi_1}\star\psi_2 = \Phi(\psi_0\star\ip{\psi_1}{\psi_2}).$$ 
(As usual when $\phi$ is an antisymmetric tricharacter, as in \cite{BHM3}, the 
invariance of the inner product renders the $\Phi$ action trivial, so that
$\Phi(\psi_0\star\ip{\psi_1}{\psi_2}) = \psi_0\star\ip{\psi_1}{\psi_2}$.)

In the associative case the dual inner products $\dip{\psi_0}{\psi_1}$ 
span the algebra $\cpt(\hilb)$ of compact operators, and in the 
nonassociative case we define the norm-closure of the span to be the twisted 
compact operators $\cpt_\phi(\hilb)$. 
The dual inner product is not generally invariant under the action of $\wh{{\sfG}}$. 
Indeed, we have 
\begin{equation}
\dip{\xi[\psi_0]}{\xi[\psi_1]}\star\xi[\psi_2] = 
\xi[\psi_0\star\ip{\psi_1}{\psi_2}] \,,
\end{equation}
from which it follows that 
$\dip{\xi[\psi_0]}{\xi[\psi_1]} = \xi\circ\dip{\psi_0}{\psi_1}\circ\xi^{-1}$.
This means that  there is a new multiplication $\star$ on $\cpt_\phi(\hilb)$, 
so that, for $K_1, K_2 \in \cpt(\hilb)$, $\psi\in \hilb$, 
\begin{equation}
(K_1\star K_2)\star\psi = \Phi(K_1\star(K_2\star\psi)) \,,
\end{equation}
where the right hand side is just the iterated natural action of 
$\cpt_\phi(\hilb)$ on $\hilb$.

\begin{lemma}
The compact operators in $\cpt_\phi(\hilb)$ are automatically adjointable with 
$\dip{\psi_2}{\psi_1}$ being the adjoint of $\dip{\psi_1}{\psi_2}$.
\end{lemma}

\begin{proof}
This is proved by consideration of 
$\dip{\psi_0}{\ip{\psi_1}{\psi_2}\star\psi_3}$ for $\psi_j \in \hilb$, 
$j = 0,1,2,3$.
To keep the orderings clearer we write $\hilb_j = \hilb$, $j=0,\ldots,3$, and 
think of $\psi_j\in \hilb_j$, so that the vectors involved in the inner 
product
$$\psi_0^*\otimes((\psi_1\otimes\psi_2^*)\otimes\psi_3)
\in \hilb_0^*\otimes((\hilb_1\otimes\hilb_2^*)\otimes\hilb_3).$$
Rewriting  $\dip{\psi_1}{\psi_2}\star\psi_3 = 
\Phi(\psi_1\star\ip{\psi_2}{\psi_3})$, is just a rebracketing, and with three 
more such steps we can rebracket it as
$(\hilb_0^*\otimes(\hilb_1\otimes\hilb_2^*))\otimes\hilb_3$,
which differs from the original by a single rebracketing, and so by just one 
application of $\Phi$ (by the pentagonal identity).
We also have
$$(\hilb_0^*\otimes(\hilb_1\otimes\hilb_2^*)) 
\cong ((\hilb_1\otimes\hilb_2^*)^*\otimes\hilb_0)^*
\cong ((\hilb_2\otimes\hilb_1^*)\otimes\hilb_0)^*,$$
from which we see that the adjoint of $\dip{\psi_1}{\psi_2}$ is 
$\dip{\psi_2}{\psi_1}$.
\end{proof}

The associative case suggests that the multiplier algebra of the twisted 
compact operators should be a twisted version of the bounded operators,
but we can be more explicit.

\begin{definition}
The twisted bounded operators $\cB_\phi(\hilb)$ are the adjointable 
operators, equipped with the $\wh{{\sfG}}$ action given by 
$\xi\cdot A = \xi\circ A\circ\xi^{-1}$, and the twisted multiplication given, 
for $A,B\in \cB_\phi(\hilb)$ and $\psi\in \hilb$, by 
\begin{equation}
(A\star B)\star\psi = \Phi(A\star(B\star\psi)) \,.
\end{equation}
\end{definition}

\begin{lemma}
If $A$ and $B$ are adjointable then so is $A\star B$ and $(A\star B)^* = 
B^*\star A^*$.
Thus the adjointable operators with this multiplication form a (generally 
nonassociative) $*$-algebra of twisted bounded operators $\cB_\phi(\hilb)$. 
The twisted compact operators $\cpt_\phi(\hilb)$  are an ideal in 
$\cB_\phi(\hilb)$ (and so, in particular, $\cpt_\phi(\hilb)$ is a subalgebra).
When $\phi$ is an antisymmetric tricharacter the twisted bounded operators 
are bounded operators but with a different multiplication.
In that case, when $\hilb = L^2({\sfG})$ with the multiplication action of $\wh{{\sfG}}$,
$\xi[\psi](x) = \xi(x)\psi(x)$, then $\cpt_\phi(\hilb)$ is the algebra 
$\cpt_\phi(L^2({\sfG}))$ defined in \cite{BHM3}.
\end{lemma}

\begin{proof}
This time we consider $\ip{\psi_1}{(A\star B)\star\psi_2}$, which arises from
$$\psi_1^*\otimes((A\otimes B)\otimes\psi_2) \in 
\hilb^*\otimes((\alg_A\otimes\alg_B)\otimes\hilb),$$
where the indices on $\alg$ just serve as reminders of where operator lives. 
Similar rebracketings to those in the previous lemma 
\begin{multline*}
\psi_1^*\otimes((A\otimes B)\otimes\psi_2)
\to \psi_1^*\otimes(A\otimes (B\otimes\psi_2))
\to (\psi_1^*\otimes A)\otimes (B\otimes\psi_2) \to\\
\to ((\psi_1^*\otimes A)\otimes B)\otimes\psi_2 
\to (\psi_1^*\otimes (A\otimes B))\otimes\psi_2 \,,
\end{multline*}
give us an element of $(\hilb^*\otimes(\alg\otimes\alg))\otimes\hilb$.
Moreover,
$$\hilb^*\otimes(\alg_A\otimes\alg_B)
\cong ((\alg_B^*\otimes\alg_A^*)\otimes\hilb)^*,$$
from which it follows that 
$$\ip{(B^*\star A^*)\star\psi_1}{\psi_2}
=\Phi(\ip{\psi_1}{(A\star B)\star\psi_2}),$$ 
so that $B^*\star A^* = (A\star B)^*$.
This proves that $A\star B$ is adjointable, and allows us to form an algebra.
In the associative case it is totally straightforward to see that 
$\cpt(\hilb)$ is an ideal because, for any $A\in \cB(\hilb)$ and $\psi_0, 
\psi_1, \psi_2 \in \hilb$, we have 
$$(A\star\dip{\psi_0}{\psi_1})\star\psi_2 = 
A\star(\psi_0\ip{\psi_1}{\psi_2}) = \dip{A\star\psi_0}{\psi_1}\star\psi_2,$$ 
showing that 
$A\star\dip{\psi_0}{\psi_1} = \dip{A\star\psi_0}{\psi_1}$. Similarly 
$\dip{\psi_0}{\psi_1} \star A= \dip{\psi_0}{A^*\star\psi_1}$, 
from which it is obvious that left and right multiplication preserve the 
generators of $\cpt(\hilb)$. 
In the nonassociative case, there are several rebracketings involved but, 
thanks to the pentagonal identity these reduce to 
$$\dip{A\star\psi_0}{\psi_1} = \Phi(A\star\dip{\psi_0}{\psi_1}).$$ 
The right hand side is 
$$\int \phi(x,y,z)\conj{\xi(x)\eta(y)\zeta(z)} 
\dip{\xi[A]\star \eta[\psi_0]}{\zeta[\psi_1]} \, dxdydz\,d\xi d\eta d\zeta,$$ 
and this is in the closure $\cpt_\phi(\hilb)$ of the span of the dual inner 
product. 

We have already seen that when $\phi$ is an antisymmetric tricharacter
the adjointable operators are the same as usual, that is, the bounded operators.
The argument of Proposition \ref{thmAd} 
applied to $\lambda(\psi_0\otimes(\psi_1\otimes\psi_2) )
= \psi_0\ip{\psi_1}{\psi_2}$ shows that $\Phi$ 
also disappears from the formula for the dual inner product in 
this case; for $\hilb = L^2({\sfG})$ with the natural inner product we get the same 
rank one operators defined by the dual inner product as when $\phi=1$, that is, 
the ordinary compact operators, which are the closure of the kernels 
$C_c({\sfG}\times {\sfG})$, and only the multiplication changes. 
More precisely we see that 
$$(\dip{\psi_0}{\psi_1}\star\psi_2)(x) 
= \Phi(\psi_0(x)\ip{\psi_1}{\psi_2}) 
 = \psi_0(x)\int_{\sfG} \phi(x,y,y)\conj{\psi_1(y)}\psi_2(y)\,dy,$$
so that for an antisymmetric tricharacter 
$\dip{\psi_0}{\psi_1}(x,y) = \psi_0(x)\conj{\psi_1(y)}$.
It is straightforward to check that multiplication follows from the formula of 
\cite[Sect 5]{BHM3}. 
\end{proof}

\begin{corollary}
The twisted bounded operators form a $*$-algebra in $\cat_{\sfG}(\phi)$ 
with $A^\star$ identified with $A^*$.
\end{corollary}

It is no coincidence that the same space $\cB(\hilb)$ can carry both 
associative and nonassociative multiplications, as one can see from the discussion of 
strictification in Appendix \ref{secappA}.

We have thus shown that the twisted bounded operators are closed under twisted multiplication, without  
deriving any  simple relationship between the norm of $A\star B$ and those of $A$ and $B$.
In particular we do not have the C$^*$-algebra identity $\| A^\star \star  A\| = \|A\|^2$.
However, there is a simple substitute for this.

\begin{proposition} \label{thmCe}
For any $A$ in the algebra of twisted bounded operators $\cB_\phi(\hilb)$, 
$A^\star\star A = 0$ if and only if $A=0$.
\end{proposition}

\begin{proof}
It is clear that $A=0$ implies $A^\star \star A = 0$.
Conversely, for any $\psi$, $\ip{A\psi}{A\psi}$ is obtained from $\ip{\psi}{(A^\star\star A)\psi}$ by the 
appropriate actions of $\Phi$. 
Since these are linear, when $A^\star\star A = 0$ we must also have 
$\ip{A\psi}{A\psi}=0$, and this forces $A=0$.
\end{proof}

\begin{definition}
A C$^*$-algebra in $\cat_{\sfG}(\phi)$ is a $*$-algebra which is $*$-isomorphic to a 
norm-closed $*$-subalgebra of $\cB_\phi(\hilb)$ for some Hilbert module $\hilb$.
(Here $*$-isomorphisms are defined as usual except that they must also be 
$\cat_{\sfG}(\phi)$-morphisms.)
\end{definition}

\section{Morita equivalence} \label{secE}

In Section \ref{secB} we defined bimodules in a category, but 
we now want to study them in a little more detail.

\begin{definition}
Let $\cX$ be a right $\alg$-module for $\alg$ a C$^*$-algebra in 
$\cat_{{\sfG}}(\phi)$, and let $\cX^*$ denote the conjugate left 
$\alg$-module. 
An  $\alg$-valued inner product on $\cX$ is a morphism 
$\cX^*\times \cX \to \alg$, written $\psi_1^*\otimes\psi_2\mapsto\rip{\psi_1}{\psi_2}$,
such that for all $\psi_1$, $\psi_2\in \cX$, and $b\in \alg$,
\begin{itemize}
\item[(i)]
$\Phi(\rip{\psi_1}{\psi_2\star b}) =\rip{\psi_1}{\psi_2}\star b$, 
where $\Phi$ is given by the action of a cocycle $\phi$ on ${\sfG}$;
\item[(ii)]
$\rip{\psi_1}{\psi_1}$ is positive in $\alg$, in the sense that it can be written as a sum of elements 
of the form $b_j^\star \star b_j$ with $b_j\in \alg$, and vanishes only when $\psi_1 =0$.
One calls $\cX$ a (right) Hilbert $\alg$-module.
\end{itemize}
\end{definition}

\begin{remark}
Assumptions (i) and  (ii) are consistent because a series of rebracketings and uses of (i) gives
$$\rip{\psi\star b}{\psi\star b} \to b^\star \star(\rip{\psi}{\psi}\star b)
= b^\star \star((\sum b_j^\star \star b_j)\star b),$$
and undoing the various rebracketings takes us back to the obviously positive term  
$$\sum_j (b_j\star b)^\star \star(b_j\star b).$$
This argument is slightly more delicate than one might realise, but as a morphism in the category 
the inner product must satisfy 
$$f[\rip{\psi_1}{\psi_2}] = \rip{f_1^*[\psi_1]}{f_2[\psi_2]},$$
as do products $f[b_j^\star\star b_j] = f_1^*[b_j^\star]\star f_2[b_j]$, so that the various actions of 
$\phi$ on the inner product and on the algebra products really do match each other.
This was the reason for taking this definition of positive rather than the other possibilities, such as 
having positive spectrum, which are equivalent in the associative case, but undefined or less useful 
here.
\end{remark}

\begin{proposition}
For all $\psi_1$ and $\psi_2\in \cX$ we have $\rip{\psi_1}{\psi_2} = \rip{\psi_2}{\psi_1}^\star$.
\end{proposition}

\begin{proof}
The polarisation identity gives
$$\rip{\psi_1}{\psi_2} = \sum_{r=0}^3 i^{-r}\ip{\psi_1 + i^r\psi_2}{\psi_1 + i^r\psi_2}
= \sum_{r=0}^3 i^{r}\ip{i^{r}\psi_1 + \psi_2}{i^{r}\psi_1 + \psi_2}$$
and, since property (ii) tells that the inner product is real, this is the same as
$\rip{\psi_2}{\psi_1}^\star$.
\end{proof}

Obviously ordinary Hilbert spaces when $\alg = \CC$, and ordinary 
Hilbert C$^*$-modules (when ${\sfG}$ is trivial), provide  examples.
Any C$^*$-algebra $\alg$, considered as an $\alg$-$\alg$-bimodule for the left 
and right multiplication actions, can be given the $\alg$-valued inner product
$\ip{a_1}{a_2}_\alg = a_1^*a_2$.
This is certainly bilinear on $\alg^*\times \alg$ and for the action by 
automorphisms one has 
$$f[\rip{a_1}{a_2}] = f[a_1^*a_2] =
f_1[a_1^*]f_2[a_2] = \rip{f_1^{*}[a_1]}{f_2[a_2]}.$$
The inner product is obviously positive, and Proposition \ref{thmCe} 
ensures that $a^\star\star a=0$  if and only if $a=0$.
This example can be combined with a Hilbert space $\hilb$ to obtain $\hilb\otimes\alg$ 
with the $\alg$-valued inner product 
$$\rip{\psi_1\otimes a_1}{\psi_2\otimes a_2} = 
\ip{\psi_1}{\psi_2}a_1^*a_2,$$
compatible with the (right) action of $\alg$ by right multiplication. 

One can similarly define $\alg$-valued inner products for left 
$\alg$-modules, either as the conjugate of the right $\alg$-module 
$\cX^*$, or directly as follows.

\begin{definition}
Let $\cX$ be a left $\alg$-module for $\alg$ a C$^*$-algebra in 
$\cat_{{\sfG}}(\phi)$, and let $\cX^*$ denote the conjugate left 
$\alg$-bimodule. 
An  $\alg$-valued inner product on $\cX$ is a morphism 
$\cX\times \cX^* \to \alg$, written $\psi_1^*\otimes\psi_2\mapsto\lip{\psi_1}{\psi_2}$,
such that for all $\psi_1$, $\psi_2\in \cX$, and $a\in \alg$,
\begin{itemize}
\item[(i)]
$\Phi_1(\lip{a\star\psi_1}{\psi_2}) = a\star\lip{\psi_1}{\psi_2}$, 
where $\Phi_1$ is given by the action of a cocycle $\phi_1$ on ${\sfG}_1$;
\item[(ii)]
$\lip{\psi_1}{\psi_1}$ is positive in $\alg$, in the sense that it can be written as a sum of elements 
of the form $a_j^\star \star a_j$ with $a_j\in \alg$, and vanishes only when $\psi_1 =0$.
One calls $\cX$ a (left) Hilbert $\alg$-module.
\end{itemize}
\end{definition}
\noindent As with the $\alg$-valued inner product a polarisation argument shows that 
$\lip{\psi_1}{\psi_2} = \lip{\psi_2}{\psi_1}^\star$.

Returning to right Hilbert C$^*$-modules, the next task is to study the algebra of adjointable 
operators $a$ on $\cX$ which commute with 
the action of $C({\sfG})$ and admit an adjoint $a^\star$ satisfying
\begin{equation}
\rip{a\star\psi}{\theta} = \Phi\rip{\psi}{a^\star \star\theta} \,.
\end{equation}
Examples are provided by rank one operators
$$\lip{\psi_0}{\psi_1}: \psi_2 \mapsto \Phi[\psi_0
\star\rip{\psi_1}{\psi_2}],$$
where $\Phi = \Phi_1\times\Phi_2$.

All this suggests an abstraction of this structure into the idea of an 
imprimitivity $\alg_1$-$\alg_2$-module  $\cX$.

\begin{definition}
Let $\cX$ be an $\alg_1$-$\alg_2$-bimodule, for $\alg_j$ a C$^*$-algebra in 
$\cat_{{\sfG}}(\phi)$, $j=1$, 2, such that the actions of $(\wh{{\sfG}},\alg_1)$ 
and $(\wh{{\sfG}},\alg_2)$are mutually commuting, and each generates the 
commutant of the other.
Let $\cX$ have $\alg_1$ and $\alg_2$-valued inner products, such that each 
algebra is adjointable in the inner product associated with the other:
\begin{equation}
\rip{a\star\psi_1}{\psi_2} =  \Phi\rip{\psi_1}{a^\star\star\psi_2},\qquad
\lip{\psi_1}{\psi_2\star b} =  \Phi\lip{\psi_1\star b^\star}{\psi_2}\,.
\end{equation}
In addition one asks that the inner products are full in the sense that their 
images are dense in $\alg_1$ and $\alg_2$, respectively, and linked by the 
imprimitivity condition that each is dual to the other,
\begin{equation}
\lip{\psi_0}{\psi_1}\star\psi_2 
= \Phi(\psi_0\star\rip{\psi_1}{\psi_2})\,.
\end{equation}
Then $\cX$ is said to be an imprimitivity bimodule for $\alg_1$ and $\alg_2$.
\end{definition}

We have defined bimodules within a single category $\cat_{\sfG}(\phi)$, but this is easily 
extended to cover
bimodules $\cX$ for C$^*$-algebras $\alg_1$ and $\alg_2$ in different tensor categories, 
$\cat_{{\sfG}_1}(\phi_1)$ and $\cat_{{\sfG}_2}(\phi_2)$, (where ${\sfG}_1$, ${\sfG}_2$ are separable 
locally compact abelian groups with three-cocycles $\phi_1$ and $\phi_2$).
We take ${\sfG} = {\sfG}_1\times {\sfG}_2$ with the product cocycle 
$\phi = \phi_1\times \phi_2$, and then $\cat_{{\sfG}_1}(\phi_1)$ forms a subcategory of 
$\cat_{\sfG}(\phi)$ on which $C({\sfG}_2)$ acts trivially (by its counit, $\epsilon_2$), 
and similarly with indices reversed.
Within $\cat_{\sfG}(\phi)$ there is no problem in taking an $\alg_1$-$\alg_2$-bimodule 
$\cX$, which as a left module is in $\cat_{{\sfG}_1}(\phi_1)$, and as a right 
module in $\cat_{{\sfG}_2}(\phi_2)$.
In particular, when ${\sfG}_2$ is trivial, one can use the bimodule to compare 
modules for twisted and untwisted algebras (cf.\ \cite{BHM3}).

As with Hilbert spaces, when $\phi_1$ and $\phi_2$ are antisymmetric 
tricharacters the rebracketing which links the dual inner products is 
independent of $\Phi$.

\begin{lemma}
In the case of $\alg_j$ in category $\cat_{{\sfG}_j}(\phi_j)$, described above, 
when $\phi_1$ and $\phi_2$ are antisymmetric tricharacters the inner products are linked by
\begin{equation}
\lip{\psi_0}{\psi_1}\star\psi_2 = \psi_0\star\rip{\psi_1}{\psi_2} \,.
\end{equation}
\end{lemma}

\begin{proof}
This follows from Proposition \ref{thmAd}, since the inner product 
$\psi_1^*\otimes\psi_2 \mapsto \rip{\psi_1}{\psi_2}$ is a $C_0({\sfG}_1)$-morphism, so that $\Phi_1$ 
acts trivially on the terms involving this, and  $\psi_1\otimes\psi_2^* \mapsto \lip{\psi_1}{\psi_2}$ 
is a $C_0({\sfG}_2)$-morphism, so that  $\Phi_2$ acts trivially on the terms involving that.
By rewriting the equation linking the inner products as
$$\Phi_2^{-1}(\lip{\psi_0}{\psi_1}\star\psi_2) 
= \Phi_1(\psi_0\star\rip{\psi_1}{\psi_2}),$$
we see that each $\Phi_j$ acts trivially, so that
$$\lip{\psi_0}{\psi_1}\star\psi_2 
= \psi_0\star\rip{\psi_1}{\psi_2}.$$ 
\end{proof}

\begin{definition}
Two C$^*$-algebras $\alg_1$ and $\alg_2$ are said to be Morita equivalent if 
there exists an imprimitivity $\alg_1$-$\alg_2$-bimodule.
\end{definition}

\begin{theorem}
Morita equivalence is an equivalence relation.
\end{theorem}

\begin{proof}
Reflexivity follows by using $\alg$ with the inner product 
$\ip{a_1}{a_2} = a_1^*a_2$ as an 
imprimitivity $\alg$-$\alg$-bimodule.
Symmetry follows by replacing an $\alg$-$\bdd$-bimodule $\cX$, by the conjugate 
$\bdd$-$\alg$-bimodule $\cX^*$, equipped with the dual inner product.

To prove transitivity, suppose that $\cX$ and $\cY$ are imprimitivity $\alg$-$\bdd$- and 
$\bdd$-$\cat$-bimodules, respectively, and set $\cZ = \cX\otimes_\bdd \cY$.
Certainly $\cZ$ is an $\alg$-$\cat$-bimodule, and it may be equipped with dual inner products
as follows.
The right inner product is obtained from the following composition of maps
\begin{multline*}
(\cY^*\otimes \cX^*)\otimes(\cX\otimes \cY)
\to \cY^*\otimes(\cX^*\otimes (\cX\otimes \cY))
\to \cY^*\otimes((\cX^*\otimes \cX)\otimes \cY)
\to \cY^*\otimes(\bdd \otimes \cY) \\
\to \cY^*\otimes \cY
\to  \cat\,,
\end{multline*}
where the first two maps are given by the appropriate associator $\Phi$, the third is the right 
inner product on $\cX$, the next is the action of $\bdd$ on $\cY$, and the last is the 
right inner product on $Y$.
The end result is $\rip{y_1}{\rip{x_1}{x_2}\star y_2}$, and this actually factors through  
$\cX\otimes_\bdd \cY$.
For example, the associator gives a map
$$\rip{x_1}{x_2\star b}\star y_2
\to  \rip{x_1}{x_2}\star (b\star y_2)$$
so that the composition depends only on $x_2\otimes_\bdd y_2$.
That will also apply for the other argument and so one really has a map $\cZ^*\otimes \cZ \to \cat$.
The map is a morphism, because
\begin{align*}
\rip{f_{(1)}^*[y_1]}{\rip{f_{(2)}^*[x_1]}{f_{(3)}[x_1]}\star f_{(4)}[y_2]}
&= \rip{f_{(1)}^*[y_1]}{f_{(2)}[\rip{x_1}{x_2}]\star f_{(3)}[y_2]}\cr
&= \rip{f_{(1)}^*[y_1]}{f_{(2)}[\rip{x_1}{x_2}\star y_2]}\cr
&= f[\rip{y_1}{\rip{x_1}{x_2}\star y_2}].
\end{align*}
The positivity follows from the positivity of the inner products on $\cX$ and $\cY$, exploiting the 
fact that exactly the same rebracketings occur for the inner products and for sums of the form 
$\sum_jb_j^\star\star b_j$.
 
The left inner product $\cZ\otimes \cZ^* \to \alg$ is similarly 
obtained from the composition of the following maps:
\begin{multline*}
(\cX\otimes \cY)\otimes(\cY^*\otimes \cX^*)
\to \cX\otimes(\cY\otimes (\cY^*\otimes \cX^*))
\to \cX\otimes((\cY\otimes \cY^*)\otimes \cX^*)
\to \cX\otimes(\bdd \otimes \cX^*) \\
\to \cX\otimes \cX^*
\to  \alg\,,
\end{multline*}
and its properties similarly checked.
\end{proof}

We saw at the end of the previous section that $\cX = \hilb$ is an imprimitivity 
$\cpt_\phi(\hilb)$-$\CC$-bimodule.

\begin{theorem}
The twisted compact operators $\cpt_\phi(\hilb)$ are Morita equivalent to 
$\CC$, with $\hilb$ providing the bimodule which gives the equivalence, 
and shows that $\cpt_\phi(\hilb)$ has trivial representation theory.
\end{theorem}

Now $\cpt_{\phi_j}(\hilb_j)$ is Morita equivalent to $\CC$ via the 
bimodule $\hilb_j$ for any Hilbert space $\hilb_j$ and twistings $\phi_j$, 
so, by the symmetry and transitivity of the equivalence, 
$\cpt_{\phi_1}(\hilb_1)$ and $\cpt_{\phi_2}(\hilb_2)$ are Morita equivalent 
via the bimodule 
$\hilb_1\otimes_\CC \hilb_2^*  = \hilb_1\otimes\hilb_2^*$.

\section{Exterior equivalence for nonassociative algebras} \label{secF}

The nonassociative algebras appear in \cite{BHM3} as 
twisted crossed products of ordinary associative algebras, 
when one had to lift an outer automorphism.
In this section we return to that situation.

In the case of a Dixmier-Douady class on a principal $\sfT$-bundle $E$ described 
by a de Rham form along the fibres, there is a homomorphism from $\sfT$ to the 
outer automorphism group of an algebra $\alg$ with spectrum $E$.
This is lifted to a map $\alpha:\sfT \to \aut(\alg)$ with 
\begin{equation}
\alpha_x\alpha_y = \ad(u(x,y))\alpha_{xy}\,,
\end{equation}
and 
\begin{equation}
u(x,y)u(xy,z) = \phi(x,y,z)\alpha_x[u(y,z)]u(x,yz)\,.
\end{equation}
Any other lifting $\beta$ would have the form
\begin{equation}
\beta_x[a] = \ad(w_x)\alpha_x[a]\,,
\end{equation}
for suitable $w_x\in \alg$, and 
\begin{equation}
\beta_x\beta_y = \ad(v(x,y))\beta_{xy}\,.
\end{equation} 
For trivial $u$ and $v$ this is just the usual exterior equivalence.

\begin{lemma}
If $\beta_x = \ad(w_x)\alpha_x$ then $\beta_x\beta_y = \ad(v(x,y))\beta_{xy}$
with 
\begin{equation}
v(x,y) = c(x,y)w_x\alpha_x[w_y]u(x,y)w_{xy}^{-1}\,,
\end{equation}
for some central $c(x,y)$, and this is the most general form of $v$.
\end{lemma}

\begin{proof}
For consistency, we must have
\begin{align*}
\ad(v(x,y)w_{xy})\alpha_{xy}[a] 
&= \ad(w_x)\alpha_x[\ad(w_y)\alpha_y[a]]\cr
&= \ad(w_x\alpha_x[w_y])[\alpha_x\alpha_y[a]]\cr
&= \ad(w_x\alpha_x[w_y]u(x,y))[\alpha_{xy}[a]],
\end{align*}
so that we must have
\begin{equation*}
v(x,y)w_{xy} = c(x,y)w_x\alpha_x[w_y]u(x,y)\,,
\end{equation*}
for some central $c(x,y)$. 
\end{proof}

For genuine representations when $u$ and $v$ are identically 1, this reduces 
to the usual requirement for exterior equivalence,
\begin{equation}
w_{xy} = w_x\alpha_x[w_y] \,.
\end{equation}

We can now prove the following result well known in the algebraic context.

\begin{lemma}
Different liftings of an outer automorphism give cohomologous cocycles $\phi$, that is, 
cocycles differing by a coboundary
\begin{equation}
(dc)(x,y,z) = \frac{c(x,y)c(xy,z)}{c(x,yz)c(y,z)} \,.
\end{equation}
\end{lemma}

\begin{proof}
In general, we can now calculate that
\begin{align*}
&v(x,y)v(xy,z) \\
&= c(x,y)c(xy,z)w_x\alpha_x[w_y]u(x,y)w_{xy}^{-1}
w_{xy}\alpha_{xy}[w_z]u(xy,z)w_{xyz}^{-1}\cr
&= c(x,y)c(xy,z)w_x\alpha_x[w_y]u(x,y)\alpha_{xy}[w_z]u(xy,z)w_{xyz}^{-1}\cr
&= c(x,y)c(xy,z)w_x\alpha_x[w_y]\alpha_x\alpha_y[w_z]u(x,y)u(xy,z)w_{xyz}^{-1}\cr
&= \phi(x,y,z)c(y,z)c(x,yz)w_x\alpha_x[w_y\alpha_y[w_z]u(y,z)]u(x,yz)
w_{xyz}^{-1}\cr
&= (\phi.dc)(x,y,z)c(y,z)c(x,yz)w_x\alpha_x[w_y\alpha_y[w_z]u(y,z)w_{yz}^{-1}]
\alpha_x[w_{yz}]u(x,yz)w_{xyz}^{-1}\cr
&= (\phi.dc)(x,y,z)w_x\alpha_x[v(y,z)]w_x^{-1}v(x,yz)\cr
&= (\phi.dc)(x,y,z)\beta_x[v(y,z)]v(x,yz),
\end{align*}
showing that $u$ and $v$ have cohomologous  cocycles $\phi$ and $\phi.dc$.
\end{proof}

\begin{theorem}
The crossed product algebras $\alg\rtimes_{\alpha,u}{\sfG}$ and $\alg\rtimes_{\beta,v}{\sfG}$, 
with $\beta_x = \ad(w_x)\alpha_x$ and $v(x,y)w_{xy} = w_x\alpha_x[w_y]u(x,y)$,
are isomorphic.
\end{theorem}

\begin{proof}
The product of $f,g\in \alg\rtimes_{\beta,v}{\sfG}$ is given by
\begin{align*}
(f\star_{\beta,v} g)(x) &= \int f(y)\beta_y[g(y^{-1}x)]v(y,y^{-1}x)\,dy\cr
&= \int f(y)w_y\alpha_y[g(y^{-1}x)]w_y^{-1}v(y,y^{-1}x)\,dy\cr
&= \int f(y)w_y\alpha_y[g(y^{-1}x)w_{y^{-1}x}]u(y,y^{-1}x)w_x^{-1}\,dy,
\end{align*}
so that we may set $f_w(x) = f(x)w_x$ to get
$$(f\star_{\beta,v} g)_w(x) = (f_w\star_{\alpha,u} g_w)(x) ,$$
and similar calculations on $f^*$ confirm that $f\mapsto f_w$ is the required isomorphism.
\end{proof}

So, up to isomorphism, the twisted crossed product depends only on the outer automorphism 
group, and a  matching multiplier.
The choice of $v$ does matter, since even when $u = 1$ the crossed product is not usually 
isomorphic to a twisted crossed product.
In fact, if we take $\phi = dc$, the composition of twisted compact operators is given by
\begin{align}
& (K_1\star K_2)(x,z) \\
&= \int K_1(x,y)K_2(y,z) dc(xy^{-1}, yz^{-1},z)\,dy\\
&= \int K_1(x,y)K_2(y,z) c(xy^{-1},yz^{-1})c(xz^{-1},z)c(xy^{-1},y)^{-1}c(yz^{-1},z)^{-1}\,dy\nonumber\\
&= c(xz^{-1},z)\int K_1(x,y)c(xy^{-1},y)^{-1}K_2(y,z)c(yz^{-1},z)^{-1} c(xy^{-1},yz^{-1})\,dy\,.\nonumber
\end{align}
This can be rewritten in terms of  $K_1^c(x,y) = K_1(x,y)c(xy^{-1},y)^{-1}$ as
\begin{equation}
(K_1\star K_2)^c(x,z) = \int K_1^c(x,y)K_2^c(y,z) c(xy^{-1},yz^{-1})\,dy\,,
\end{equation}
showing that the new multiplication is still twisted.
Indeed, it is obtained by applying the untwisted multiplication to the image of $K_1\otimes K_2$
under the action of $c\in C({\sfG})\otimes C({\sfG})$.
Even such slightly deformed products can have very different differential calculi, as investigated 
by Majid and collaborators (particularly the recent preprint \cite{BM4}).

We now define two automorphism groups with algebra-valued cocycles 
$(\alpha,u)$ and $(\beta,v)$ to be exterior equivalent if 
\begin{equation}
\beta_x = \ad(w_x)\alpha_x, \qquad\hbox{{\rm and}}\qquad
v(x,y) = w_x\alpha_x(w_y)u(x,y)w_{xy}^{-1}\,.
\end{equation}
These can be rewritten as
\begin{equation}
\alpha_x = \ad(w_x^{-1})\beta _x, \qquad\hbox{{\rm and}}\qquad
u(x,y) = \beta_x(w_y^{-1})w_x^{-1}v(x,y)w_{xy}\,,
\end{equation}
demonstrating the reflexivity of the equivalence, and similarly the product of 
the algebra elements gives transitivity.

We now can easily rephrase (and shorten) the derivation of Packer--Raeburn 
equivalence \cite[Cor 4.2]{BHM3}.
The cocycle equation for $u$ can be written as
$$u(x,y)u(xy,y^{-1}x^{-1}z) 
= \phi(x,y,y^{-1}x^{-1}z)\alpha_x[u(y,y^{-1}x^{-1}z)]u(x,x^{-1}z).$$
We now work in $\alg\otimes\cpt(L^2({\sfG}))$ which acts on 
$a\in \alg\otimes L^2({\sfG})$.
We define 
$$(w_xa)(z) = u(x,x^{-1}z)a(x^{-1}z),$$
and $\beta_x = \ad(w_x)\alpha_x$
and use the above version of the cocycle identity to show that 
$(\alpha,u)$ is exterior equivalent to $(\beta,v)$ with 
$(v(x,y)a)(z) = \phi(x,y,y^{-1}x^{-1}z)a(z)$.

\section{The general duality result} \label{secG}

In this section we generalise the construction at the 
end of Section \ref{secF}  to general C$^*$-algebras with twisted 
actions. 

\begin{theorem}\label{thm:takai}
Let $\bdd$ be a C$^*$-algebra on which the group ${\sfG}$ acts by twisted automorphisms 
$\beta_g$, with twisting given by $v(x,y)$, satisfying the deformed cocycle 
condition with tricharacter obstruction $\phi$. The twisted crossed product 
$(\bdd\rtimes_{\beta,v}{\sfG})\rtimes \wh{{\sfG}}$ is isomorphic to the algebra of $\bdd$-
valued twisted kernels $\bdd\otimes \cpt_{\conj{\phi}}(L^2({\sfG}))$ with the product 
\begin{equation}
(k_1\star k_2)(w,z) 
= \int_{{\sfG}} k_1(w,u)k_2(u,z)\phi(wu^{-1},uz^{-1},z)^{-1}\,du \,.
\end{equation}
The double dual action on $(\bdd\rtimes_{\beta,v}{\sfG})\rtimes \wh{{\sfG}}$ is equivalent to 
\begin{equation}
(\wh{\wh{\beta}}_yk_F)(w,z)
=\phi(wz^{-1},z,y)\ad(V_y)^{-1}[\beta_{y}[k_F(wy,zy)]]\,,
\end{equation}
which is the product of the original action $\beta_y$ and a twisted adjoint action of
\begin{equation}
V_y(z) = \frac{\phi(y,z^{-1},z)}{\phi(yz^{-1},zy^{-1},y)}v(y,z^{-1})\,,
\end{equation}
on $\bdd$-valued kernels. 
\end{theorem}

\begin{proof}
The twisted crossed product $\bdd\rtimes_{\beta,v} {\sfG}$ consists of 
$\bdd$-valued functions on ${\sfG}$ with product
$$(f*g)(x) = \int_{\sfG} f(y)\beta_y[g(y^{-1}x)]v(y,y^{-1}x)\,dy,$$
and $(\bdd\rtimes_{\beta,v}{\sfG})\rtimes \wh{{\sfG}}$ consists of $\bdd$-valued functions 
on ${\sfG}\times\wh{{\sfG}}$ with product 
\begin{align*}
(F*G)(x,\xi) &= \int_{\wh{G}} 
(F(.,\eta)\wh{\beta}_\eta[G(.,\eta^{-1}\xi)])(x)\,d\eta\cr &= 
\int_{\wh{G}\times G} F(y,\eta)\beta_y\wh{\beta}_\eta [G(y^{-1}x,\eta^{-
1}\xi)]v(y,y^{-1}x)\,dyd\eta\cr &= \int_{\wh{G}\times G} 
F(y,\eta)\beta_y[\eta(y^{-1}x) G(y^{-1}x,\eta^{-1}\xi)])v(y,y^{-
1}x)\,dyd\eta.
\end{align*}
We now Fourier transform with respect to the second 
argument, so that $$\wh{F}(x,z) = \int_{\wh{G}} F(x,\xi)\xi(z)\,d\xi,$$ to get 
\begin{align*}
(\wh{F*G})(x,z) 
&= \int F(y,\eta)\beta_y[\eta(y^{-1}x)
G(y^{-1}x,\eta^{-1}\xi)])v(y,y^{-1}x)\xi(z)\,dyd\eta d\xi\cr
&= \int F(y,\eta)\eta(y^{-1}xz)\beta_y[G(y^{-1}x,\eta^{-1}\xi)
(\eta^{-1}\xi)(z))]v(y,y^{-1}x)\,dyd\eta d\xi\cr
&= \int_{G} \wh{F}(y,y^{-1}xz)\beta_y[\wh{G}(y^{-1}x,z)]v(y,y^{-1}x)\,dy.
\end{align*}
Next we introduce 
$\wh{k}_F(w,z) = \beta_{w^{-1}}[\wh{F}(wz^{-1},z)]v(w^{-1},wz^{-1})$
and by setting $w = xz$ in the last product formula, applying $\beta_{w^{-
1}}$, and using the standard identities for $v$ and $\phi$, we obtain 
\begin{align*}
&\wh{k}_{F*G}(w,z) 
= \beta_{w^{-1}}[\wh{F*G})(wz^{-1},z)]v(w^{-1},wz^{-1})\cr 
&= \int_{G} \beta_{w^{-1}}[\wh{F}(y,y^{-1}w)] \beta_{w^{-1}}
\beta_y[\wh{G}(y^{-1}wz^{-1},z)] \beta_{w^{-1}}[v(y,y^{-1}wz^{-1})]
v(w^{-1},wz^{-1})\,dy\cr &= \int_{G} \beta_{w^{-1}}[\wh{F}(y,y^{-1}w)] 
v(w^{-1},y)\beta_{w^{-1}y}[\wh{G}(y^{-1}wz^{-1},z)] \cr
& \qquad\times v(w^{-1}y,y^{-1}wz^{-1})
\phi(w^{-1},y,y^{-1}wz^{-1})\,dy\cr 
&= \int_{G} \wh{k}_F(w,y^{-1}w)\wh{k}_G(y^{-1}w,z)
\phi(w^{-1},y,y^{-1}wz^{-1})\,dy\cr
&= \int_{G} \wh{k}_F(w,u)\wh{k}_G(u,z)\phi(w^{-1},wu^{-1},uz^{-1})\,du.
\end{align*}
Now, by the cocycle identity 
$$\phi(w^{-1},wu^{-1},uz^{-1})
= \frac{\phi(u^{-1},uz^{-1},z)\phi(w^{-1},wu^{-1},u)}
{\phi(wu^{-1},uz^{-1},z)\phi(w^{-1},wz^{-1},z)},$$
from which it follows that with 
$k_F(w,u) = \wh{k}_F(w,u)\phi(w^{-1},wu^{-1},u)$ we have
$$k_{F*G}(w,z) 
= \int_{G} k_F(w,u)k_G(u,z)\phi(wu^{-1},uz^{-1},z)^{-1}\,du.$$
Thus we have an isomorphism with $\bdd$-valued twisted kernels 
$\bdd\otimes \cpt_{\conj{\phi}}(L^2(\sfG))$.

We can also compute that the double dual (untwisted) action of $y\in G$ takes 
$F(x,\xi)$ to $\xi(y)F(x,\xi)$. Integrating against $\xi(z)$ we see that this 
takes $\wh{F}(x,z)$ to $\wh{F}(x,zy)$, and so $\wh{k}_F$ to 
\begin{align*}
&(\wh{\wh{\beta}}_y\wh{k}_F)(w,z) = \beta_{w^{-1}}[\wh{F}(wz^{-1},zy)]v(w^{-1},wz^{-1})\cr
&= \beta_{y(wy)^{-1})}[\wh{F}((wy)(zy)^{-1},zy)]
v(w^{-1},wz^{-1})\cr
&= \ad(v(y,(wy)^{-1}))^{-1}\beta_{y}
\beta_{(wy)^{-1}}[\wh{F}((wy)(zy)^{-1},zy)]v(w^{-1},wz^{-1})\cr
&= v(y,(wy)^{-1})^{-1}\beta_y[\wh{k}_F(wy,zy)
v((wy)^{-1},wz^{-1})^{-1}]v(y,(wy)^{-1})v(w^{-1},wz^{-1})\cr
&= v(y,(wy)^{-1})^{-1}\beta_{y}[\wh{k}_F(wy,zy)]
\beta_{y}[v((wy)^{-1},wz^{-1})]^{-1}v(y,(wy)^{-1})v(w^{-1},wz^{-1})\cr
&= v(y,(wy)^{-1})^{-1}\beta_{y}[\wh{k}_F(wy,zy)]
\phi(y,(wy)^{-1},wz^{-1})^{-1}v(y,(zy)^{-1}).
\end{align*}

This in turn, with a couple of applications of the pentagonal identity, gives 
the following expression for $(\wh{\wh{\beta}}_yk_F)(w,z)$:
\begin{align*}
&\phi(w^{-1},wz^{-1},z)\phi(y^{-1}w^{-1},wz^{-1},zy)^{-1}
\phi(y,(wy)^{-1},wz^{-1})^{-1}v(y,(wy)^{-1})^{-1}\\
& \qquad \times \beta_{y}[k_F(wy,zy)]v(y,(zy)^{-1})\cr
&=\frac{\phi(w^{-1},wz^{-1},z)}{\phi(w^{-1},wz^{-1},zy)}
\frac{\phi(y,y^{-1}z^{-1},zy)}{\phi(y,y^{-1}w^{-1},wy)}
v(y,(wy)^{-1})^{-1}\beta_{y}[k_F(wy,zy)]v(y,(zy)^{-1})\cr
&=\phi(wz^{-1},z,y)\frac{\phi(w^{-1},w,y)}{\phi(y,y^{-1}w^{-1},wy)}
v(y,(wy)^{-1})^{-1}\beta_{y}[k_F(wy,zy)]v(y,(zy)^{-1}) \cr
& \qquad \times \frac{\phi(y,y^{-1}z^{-1},zy)}{\phi(z^{-1},z,y)}
=\phi(wz^{-1},z,y)\ad(V_y)^{-1}\beta_{y}[k_F(wy,zy)],
\end{align*}
where $V_y$ is defined in the Theorem.
This is the product of the original action $\beta_y$, and a twisted action on 
$\bdd$-valued kernels, which combines the original action $\beta_y$ with an 
adjoint action of $V_y$, and an action on kernels of the type discussed in 
\cite[Sect 5]{BHM3}. 
(Two actions $\alpha$ and $\beta$ linked by an inner automorphism $\ad\, u$ 
can be removed by an exterior equivalence.)
\end{proof}

As a consequence of this result we can, when convenient, replace an algebra 
$\bdd$ with twisted action by a (stable) nonassociative algebra $\bdd\otimes 
\cpt_\phi(L^2({\sfG}))$ with an ordinary action. 

Since the arguments are now quite general they could also be used to show that 
the third dual is isomorphic to the first dual tensored with ordinary compact 
operators.

\section{Twisted and repeated crossed products} \label{secH}

The standard Takai duality  theorem can be seen from a different perspective 
by identifying the repeated crossed product 
$(\alg\rtimes {\sfG})\rtimes \wh{{\sfG}}$ with the twisted crossed product
$\alg \rtimes({\sfG}\times_\sigma\wh{{\sfG}})$, 
where the twisting is done by the Mackey multiplier $\sigma$ on 
${\sfG}\times \wh{{\sfG}}$ given by $\sigma((y,\eta),(x,\xi)) = \eta(x)$.
(For ${\sfG} = \RR$,  every multiplier is equivalent to a Mackey multiplier.)
In fact the crossed product $\alg\rtimes {\sfG}$ consists of 
$\alg$-valued functions $a$, $b$, on ${\sfG}$ with product 
\begin{equation}
(a*b)(x) = \int_{\sfG} a(y)\alpha_y[b(y^{-1}x)]\,dy\,,
\end{equation}
with $\wh{{\sfG}}$-action $(\wh{\alpha}_\xi[a])(x) = \conj{\xi(x)}a(x)$,
and the repeated crossed product similarly consists of functions from 
${\sfG}\times_\sigma\wh{{\sfG}}$ to $\alg$ with product
\begin{equation}
(a*b)(x,\xi) = \int a(y,\eta)\conj{\eta(y^{-1}x)}\alpha_y[b(y^{-1}x,\eta^{-1}\xi)]\,dyd\eta \,.
\end{equation}
Now a twisted crossed product given by a projective action $\beta$ of 
${\sfG}\times_\sigma\wh{{\sfG}}$ 
would have product
\begin{equation}
(a*b)(x,\xi) = \int a(y,\eta)\beta_{(y,\eta)}[b(y^{-1}x,\eta^{-1}\xi)]
\conj{\sigma((y,\eta),(y^{-1}x,\eta^{-1}\xi))}\,dyd\eta\,,
\end{equation}
and we see that these match if $\beta_{(y,\eta)} = \alpha_y$ and 
$\sigma((y,\eta),(x,\xi)) = \eta(x)$.
The duality theorem then follows from Green's results on imprimitivity algebras
and the fact that the twisted group algebra of ${\sfG}\times_\sigma\wh{\sfG}$ 
is essentially the algebra of compact operators on $L^2({\sfG})$.

This equivalence between a twisted crossed product and a repeated (untwisted) 
crossed product has an analogue when the twisting is done with a three-cocycle,
which can be exploited in reverse, to reinterpret the twisted crossed product 
as a repeated crossed product.
We take two exterior equivalent twisted automorphism groups 
$(\beta,v)$ and $(\wh{\beta},\wh{v})$ of a C$^*$-algebra $\alg$:
\begin{equation}
\wh{\beta}_x = \ad(w_x)\beta_x, \qquad\hbox{{\rm and}}\qquad
\wh{v}(x,y) = w_x\beta_x(w_y)v(x,y)w_{xy}^{-1}\,.
\end{equation}
Both $v$ and $\wh{v}$ define the same three-cocycle $\phi$, which we 
shall assume to be an antisymmetric tricharacter, so that when trivial it is 
identically 1.

\begin{theorem}
For a separable locally compact abelian group $\sfG$ and a stable C$^*$-algebra
$\alg$, let $\beta: \sfG \to \aut(\alg)$ be a twisted homomorphism with 
$\beta_\xi\beta_\eta = \ad(v(\xi,\eta))\beta_{\xi\eta}$, for all 
$\xi,\eta\in \sfG$.     
Suppose that $\sfG = \sfG_1\times \sfG_2$, where $\sfG_j$ has trivial Moore 
cohomology $H^3(\sfG_j,\mathsf U(1))$.
Suppose that the corresponding three-cocycle $\phi$ is identically 1 on 
$\sfG_1$ and is also 1 when two of its arguments are in $\sfG_2$
(this being true for antisymmetric tricharacters on $\RR^2\times\RR$). 
Then we may take $v$ to be trivial on the subgroups $\sfG_j$ and 
$(\beta,v)$ is exterior equivalent to $(\wh{\beta},\wh{v})$, where 
\begin{equation}
\wh{\beta}_{xX} = \beta_x\beta_X, \qquad {{\rm and}}\qquad
\wh{v}(xX,yY) = \beta_x[\wt{v}(X,y)]
\end{equation}
where lower case letters denote elements of $\sfG_1$ and capitals elements of 
$\sfG_2$, and $\wt{v}(X,y) = v(X,y)v(y,X)^{-1}$. 
The cocycle $\wt{v}$ satisfies
\begin{align}
\wt{v}(XY,z) & = \beta_X[\wt{v}(Y,z)]\wt{v}(X,z),\qquad {{\rm and}}\qquad \nonumber \\
\wt{v}(X,yz)  & = \phi(X,y,z)^{-3}\wt{v}(X,y)\beta_y[\wt{v}(X,z)]\,,
\end{align}
and $\wh{v}$ defines the three-cocycle $\varphi(xX,yY,zZ) = \phi(X,y,z)^3$.
\end{theorem}

\begin{proof}
Since $\phi= 1$ on the subgroups, the stabilisation theorem 
(\cite[Cor4.2]{BHM3} and below) tells us that we may also take $v$ to be 1 
on the subgroups.
It follows from the definitions that 
$\wh{\beta}_{xX} = \ad(v(x,X))\beta_{xX}$, so that we could take $w_{xX} = v(x,X)$.
In fact it is more convenient to make a slightly different choice since we also have
$$\beta_\xi\beta_\eta = v(\xi,\eta)\beta_{\xi\eta} 
= \ad(\wt{v}(\xi,\eta))\beta_\eta\beta_\xi,$$
from which we deduce
\begin{align*}
(\beta_x\beta_X)(\beta_y\beta_Y) 
&= \beta_x(\beta_X\beta_y)\beta_Y\\
&= \beta_x\ad(\wt{v}(X,y))\beta_y\beta_X\beta_Y\\
&= \ad(\beta_x[\wt{v}(X,y)])\ad(v(x,y))\beta_{xy}\ad(v(X,Y))\beta_{XY}\,.
\end{align*}
Since the cocycle $\wh{v}$ is determined up to scalars, and $v(x,y) = 1$, 
$v(X,Y) = 1$ on the subgroups, this gives us the required form
$\wh{v}(xX,yY) = \beta_x[\wt{v}(X,y)]$.
Alternatively we can see that up to scalar factors
$w_{xX}\beta_{xX}[w_{yY}]v(xX,yY)w_{xyXY}^{-1}$ is
$$\beta_x[\wt{v}(X,y)]v(x,y)\beta_{xy}][v(X,Y)]
= \beta_x[\wt{v}(X,y)]$$
using the triviality of the restriction of $v$ to the subgroups.
Now, we also have
\begin{align*}
& v(XY,z) = \phi(X,Y,z)v(X,Y)^{-1}\beta_X[v(Y,z)]v(X,Yz)\\
&= \phi(X,Y,z)v(X,Y)^{-1}\beta_X[\wt{v}(Y,z)]\phi(X,z,Y)^{-1}v(X,z)v(Xz,Y)\\
&= \phi(X,Y,z)\phi(X,z,Y)^{-1}v(X,Y)^{-1}\\
& \qquad \times \beta_X[\wt{v}(Y,z)]\wt{v}(X,z)
\phi(z,X,Y)\beta_z[v(X,Y)]v(z,XY)\,,
\end{align*}
and using the triviality of $v$ on subgroups, and the fact that $\phi$ must be 
trivial whenever two of its arguments are in $\sfG_2$, this reduces to
$$\wt{v}(XY,z)= \beta_X[\wt{v}(Y,z)]\wt{v}(X,z).$$
Similarly, we have
\begin{align*}
&\beta_X[v(y,z)]v(X,yz) 
= \phi(X,y,z)^{-1}v(X,y)v(Xy,z)\\
&= \phi(X,y,z)^{-1}\phi(y,X,z)\wt{v}(X,y)\beta_y[v(X,z)]v(y,Xz)\\
&= \phi(X,y,z)^{-1}\phi(y,X,z)\wt{v}(X,y)\beta_y[\wt{v}(X,z)]
\phi(y,z,X)^{-1}v(y,z)v(yz,X)\,,
\end{align*}
and using the triviality of $v$ on subgroups as well as the antisymmetry of 
$\phi$ this reduces to 
$$\wt{v}(X,yz) 
= \phi(X,y,z)^{-3}\wt{v}(X,y)\beta_y[\wt{v}(X,z)].$$
We then have
\begin{align*}
&\wh{v}(xX,yY)\wh{v}(xyXY,zZ)\wh{v}(xX,yzYZ)^{-1} \\
&= \beta_x[\wt{v}(X,y)]\beta_{xy}[\wt{v}(XY,z)]\beta_x[\wt{v}(X,yz)]^{-1}\\
&= \beta_x[\wt{v}(X,y)\beta_y[\wt{v}(XY,z)]\wt{v}(X,yz)^{-1}]\\
&= \phi(X,y,z)^{-3}\beta_x[\wt{v}(X,y)\beta_y[\wt{v}(XY,z)\wt{v}(X,z)^{-1}
\wt{v}(X,y)^{-1}]\\
&= \phi(X,y,z)^{-3}\beta_x[\wt{v}(X,y)\beta_y[\beta_X[\wt{v}(Y,z)]
\wt{v}(X,y)^{-1}]\\
&= \phi(X,y,z)^{-3}\beta_x\beta_X[\beta_y[\wt{v}(Y,z)]]\\
&= \phi(X,y,z)^{-3}\beta_x\beta_X[\wh{v}(yY,zZ)]\,,
\end{align*}
showing that a cocycle identity holds for $\wh{v}$.
The corresponding three-cocycle $\varphi(xX,yY,zZ) = \phi(X,y,z)^3$ is not 
antisymmetric, but its antisymmetrisation,
\begin{align*}
[\varphi(xX,yY,zZ)\varphi(yY,zZ,xX)\varphi(zZ,xX,yY)]^{\frac13} &=
\phi(X,y,z)\phi(Y,z,x)\phi(Z,x,y)\\
&= \phi(xX,yY,zZ),
\end{align*}
is just the original cocycle $\phi$.
\end{proof}

\begin{theorem}
Writing $f_Y(y) = f(yY)$, the crossed product algebra 
$\alg\rtimes_{\wh{\beta},\wh{v}}G \sim (\alg\rtimes {\sfG}_1)\rtimes_{\wt{\beta}} {\sfG}_2$ 
has twisted convolution product
\begin{equation}
(f*g)_X = \int_{{\sfG}_2} f_Y *_{{}_1}\wt{\beta}_Y[g_{Y^{-1}X}]\,dY\,,
\end{equation}
where $\wt{\beta}_X[f](y) = \beta_X[f(y)]\wt{v}(X,y)$, and $*_{{}_2}$ denotes the 
convolution product on $\alg\rtimes {\sfG}_1$.
The map $X \to \wt{\beta}_X$ is a group homomorphism, and $\wt{\beta}_X$ gives an  
isomorphism between twisted crossed products for 
multipliers $\sigma$ and $\phi_X\sigma$, where $\phi_X(y,z) = \phi(X,y,z)$.
\end{theorem}

\begin{proof}
We calculate that 
\begin{align*}
(f*g)(xX) 
&= \int f(yY)\beta_y\beta_Y[g(y^{-1}xY^{-1}X)
\beta_y[\wt{v}(Y,y^{-1}x)]\,dydY\cr
&= \int f(yY)\beta_y[\beta_Y[g(y^{-1}xY^{-1}X)]\wt{v}(Y,y^{-1}x)]\,dydY\\
&= \int f(yY)\beta_y[\wt{\beta}_Y[g(y^{-1}xY^{-1}X)]]\,dydY\,,
\end{align*}
so that
$$(f*g)_X = \int f_Y *_{{}_1}\wt{\beta}_Y[g_{Y^{-1}X}]\,dY.$$
We can then check that
\begin{align*}
\wt{\beta}_X\wt{\beta}_Y[f](z) 
&= \beta_X[\beta_Y[f(z)]\wt{v}(Y,z)]\wt{v}(X,z)\\
&= \beta_{XY}[f(z)]\beta_X[\wt{v}(Y,z)]\wt{v}(X,z)\\
&= \beta_{XY}[f(z)]\wt{v}(XY,z)\\
&= \wt{\beta}_{XY}[f](z)\,.
\end{align*}
Next consider a twisted crossed product with $U(1)$-valued multiplier $\sigma$
\begin{align*}
(\wt{\beta}_X[f]*\wt{\beta}_X[g])(x)
&= \int \wt{\beta}_X[f](y)\beta_y[\wt{\beta}_X[g](y^{-1}z)]
\sigma(y,y^{-1}x)\,dy\cr
&= \int \beta_X[f(y)]\wt{v}(X,y)\beta_y[\beta_X[g(y^{-1}z)]\wt{v}(X,y^{-
1}z)]\sigma(y,y^{-1}x)\,dy\cr
&= \int \beta_X[f(y)]\beta_X[\beta_y[g(y^{-1}x)]]
\wt{v}(X,y)\beta_y[\wt{v}(X,y^{-1}x)]\sigma(y,y^{-1}x)\,dy\cr
&= \beta_X[\int f(y)\beta_y[g(y^{-1}x)]\varphi(X,y,y^{-1}x)\sigma(y,y^{-1}x)
\,dy]\wt{v}(X,x)\,,
\end{align*}
which differs from $\wt{\beta}_X[f*g](x)$ by the insertion of 
$\varphi(X,y,y^{-1}x)$ in the convolution integral, changing the multiplier 
$\sigma$ to $\varphi(X,\cdot,\cdot)\sigma$.
\end{proof}

This result tells us that the twisted crossed product multiplication for 
$\alg\rtimes \sfG$ can be obtained by doing repeated crossed products but with 
a modified automorphism action of the final subgroup.
The result is still nonassociative because one still has the three-cocycle 
$\varphi$. 
There is no inconsistency because $\wt{\beta}_X$ does not 
act as automorphisms of a crossed product $\alg\rtimes \sfG_1$.
We shall now show how to modify things to get a more useful result.

Suppose now that we replace the algebra $\alg$ by an algebra $\cb$ which 
admits not only an action of ${\sfG}$ but also a compatible action of 
$C(\sfG_2) = C(\sfG/\sfG_1)$. (This would be automatic for an algebra induced to $\sfG$ 
from a subgroup ${\sfN}$.)

\begin{theorem}
Let $\cb$ be a C$^*$-algebra admitting an action $\beta$ of the abelian group 
$\sfG = \sfG_1\times \sfG_2$ with multiplier $v$ and $\phi$ satisfying the conditions of 
the previous theorems and also a compatible action of $C({\sfG}_2)$.
Then $\cb\rtimes_{\beta,v} \sfG$ is stably equivalent to an ordinary repeated 
crossed product $(\cb\rtimes_{\wh{\beta}_1}{\sfG}_1)\rtimes_{\wh{\beta}_2} {\sfG}_2$, 
where $\wh{\beta}_j$ is the restriction of $\wh{\beta}$ to ${\sfG}_j$.
\end{theorem}

\begin{proof}
We have already seen that, by replacing $(\beta,v)$ by $(\wh{\beta},\wh{v})$,
the multipliers on the subgroups are trivial and that when one splits the 
crossed product into a repeated crossed product the first part is just an 
ordinary crossed product with the action $\wh{\beta}_1$.
Now, by the Packer--Raeburn trick, this is stably equivalent to a projective 
crossed product with the $C({\sfG}_2)$-valued multiplier $\varphi_\bullet(y,z)$ 
defined as the function $X\mapsto \varphi(X,y,z)$,
The advantage is that $\wt{\beta}_X$ is now an automorphism, since when we put 
$\sigma = \varphi_\bullet$ we have
\begin{align*}
& (\wt{\beta}_X[f]*_{\varphi_\bullet}\wt{\beta}_X[g])(x)
= \int \wt{\beta}_X[f](y)\beta_y[\wt{\beta}_X[g](y^{-1}z)]
\varphi_\bullet(y,y^{-1}x)\,dy\\
&= \int \beta_X[f(y)]\wt{v}(X,y)\beta_y[\beta_X[g(y^{-1}z)]
\wt{v}(X,y^{-1}z)]\varphi_\bullet(y,y^{-1}x)\,dy\\
&= \int \beta_X[f(y)]\beta_X[\beta_y[g(y^{-1}x)]]
\wt{v}(X,y)\beta_y[\wt{v}(X,y^{-1}x)]\varphi_\bullet(y,y^{-1}x)\,dy\\
&= \int \beta_X[f(y)]\beta_X[\beta_y[g(y^{-1}x)]]\varphi_\bullet(y,y^{-1}x)
\varphi(X,y,y^{-1}x) \,dy]\wt{v}(X,x)\\
&= \beta_X[\int f(y)\beta_y[g(y^{-1}x)]\varphi_\bullet(y,y^{-1}x)\,dy]
\wt{v}(X,x) = \wt{\beta}_X[f*_{\varphi_\bullet}g]\,,
\end{align*}
where 
$\varphi_Z(y,z)\varphi(X,y,z) = \varphi_{ZX}(y,z) 
= \wh{\beta}_X(\varphi_Z(y,z)\wh{\beta}_X^{-1}$
follows from compatibility of the actions. 
\end{proof}

This result takes us quite a long way towards proving the analogue of the 
Connes-Thom isomorphism in our context.
We take ${\sfG} = \RR^3$, with the subgroups ${\sfG}_1 = \RR^2$, and 
${\sfG}_2 = \RR$.
In \cite{BHM3} the algebra has the form $\cb = \ind_{\integer^3}^{\RR^3}\alg$, 
and so has an action of $\RR^3$.
As already noted the Moore cohomology group $H^3$ is trivial on the subgroups 
${\sfG}_1 = \RR^2$ and ${\sfG}_2 = \RR$.
Using Theorem 7.1 we express the twisted crossed product 
$\cb\rtimes\RR^3$ as a repeated crossed product 
$(\cb\rtimes\RR^2)\rtimes\RR$.
The standard Connes-Thom isomorphism tells us that 
$K_*(\cb\rtimes\RR^2) \cong K_*(\cb)$, so that the first crossed product 
does not change the $K$-theory.
The second crossed product with ${\sfG}_2$ is more problematic because the 
group does not act as automorphisms.
However, it is stably equivalent to the case where one does have automorphisms.
In that stably equivalent algebra the ordinary Connes-Thom theorem then 
asserts that the $K$-theory of the crossed product 
$(\cb\rtimes\RR^2)\rtimes\RR$ is the same as that of 
$\cb\rtimes\RR^2$ and so of $\cb$, apart from a shift of 1 
in degree.
Superficially this appears to have proved the desired generalisation of the 
Connes-Thom theorem to the twisted algebra $\cb\rtimes_{\beta,v} \RR^3$, 
but the missing ingredient is to check that the stabilisation equivalence 
valid for ordinary $K$-theory is also consistent with the definition of 
$K$-theory in the category $\cat_{\sfG}(\phi)$.

\section{Some remarks and speculations related to $K$-theory} \label{secI}

Here we give an application, to a version of $K$-theory, of Takai duality in our context.
For $\alg$ an algebra in the category $\cat_{\sfG}(\phi)$ we can define the ring 
$K_0(\alg,\Phi)$ of stable equivalence classes of projective finite rank 
$\alg$-modules in $\cat_{\sfG}(\phi)$.
We have already seen that these modules are also $\alg\rtimes \wh{{\sfG}}$-modules,
and it follows from Proposition A.2 that these are also finite rank and 
projective. 
There is thus a natural identification of $K_0(\alg,\Phi)$ with the stable 
equivalences of finite rank projective $\alg\rtimes \wh{{\sfG}}$-modules.
There is a caveat that $K_{0}(\alg,1)$ does not reduce to the $K$-theory of the $\sfG$-algebra $\alg$,
since the finite projective $\alg$-modules used in the definition
are also expected to be $\sfG$-invariant. The reason is that only $\sfG$-invariant
projections are well defined in the category $\cat_{\sfG}(\phi)$, cf. below.
In fact, $K_0(\alg,1) \cong K_0(\alg\rtimes \wh{{\sfG}})$, and for ${\sfG} \cong \RR^d$ 
Connes-Thom isomorphism theorem gives 
$K_0(\alg\rtimes\wh{{\sfG}}) \cong K_{d}(\alg)$, so that 
$K_0(\alg,1) \cong K_{d}(\alg)$.

We actually want to apply this to a C$^*$-algebra of the form $\alg = \bdd\rtimes {\sfG}$, 
with $\bdd$ associative, to find $K_0(\bdd\rtimes {\sfG},\Phi)$ in terms of 
equivalence classes of $(\bdd\rtimes {\sfG})\rtimes\wh{{\sfG}}$-modules.
By Theorem \ref{thm:takai}, one has 
$(\bdd\rtimes {\sfG})\rtimes\wh{{\sfG}} \cong \bdd\otimes\cpt_{\conj{\phi}}(L^2({\sfG}))$, 
which by Appendix A, can be strictified to the associative C$^*$-algebra 
$\bdd\otimes\cpt(L^2({\sfG}))$.
(There is a natural correspondence between the modules since the 
nonassociative effect of $\Phi$ appears only for repeated actions, and the 
action itself can be defined in the same way for both cases.
The strictification functor thus preserves the properties of being finite rank 
and projective.)
This associative C$^*$-algebra is Morita equivalent to $\bdd$ itself, so that the
stable equivalence classes of finite rank projective modules for 
$(\bdd\rtimes {\sfG})\rtimes\wh{{\sfG}}$ are in natural bijective correspondence with those for 
$\bdd$, that is, with $K_0(\bdd\otimes\cpt(L^2({\sfG})))$.
In other words 
$K_0(\bdd\rtimes {\sfG},\Phi) \cong K_0(\bdd\otimes\cpt(L^2({\sfG})))$.

At first sight a degree change appears to be missing, but this is a 
consequence of the way we have defined $K_0(\alg, \Phi)$, as previously explained.
The morphisms in the category $\cat_\sfG(\phi)$ with associator $\phi$ are 
${\sfG}$-maps, and so projective modules for an algebra $\alg$ in this category are 
submodules of free modules defined by ${\sfG}$-invariant projections $E$.
Such modules can also be thought of as submodules of free modules defined by 
idempotents $e$ in a matrix algebra $M_n(\alg)$, so that 
$Ev = v\star e$.
We then see that
\begin{equation}
v\star g[e] = g(g^{-1}v\star e) = g(Eg^{-1}v) = Ev = v\star e \,,
\end{equation}
so that $e = g[e]$ is invariant under the action of ${\sfG}$.

To properly define the $K$-theory
of a C$^*$-algebra in the category $\cat_\sfG(\phi)$, we need to consider 
an enveloping category that includes not just $\sfG$-morphisms. One such candidate 
is the Karoubian enveloping category, and will be considered in 
a future work. We also plan to investigate Tannakian duality and its 
consequences in our context. More precisely, let $\sfG$ 
denote the Euclidean group, and $\phi$ the 3-cocycle  on it
as in the text. Then our
tensor category  $\cat_\sfG(\phi)$ is just a twisted
representation category,  ${\rm Rep}_c(\sfG, \phi)$, where the subscript
 is a reminder that we take topology into
consideration. Then the dual tensor category
consists of continuous, tensorial functors
$
\cF : {\rm Rep}_c(\sfG, \phi) \to \cB(V),
$
where $\cB(V)$ denotes bounded operators on the Hilbert space $V$.
Here $V$ varies, and tensorial means respecting the structures
in the tensor category. The dual tensor category is a tensor category
denoted by ${\rm Rep}_c(\sfG, \phi)'$, and the putative analog of
Tannakian duality in this context would say that
$
 {\rm Rep}_c(\sfG, \phi)$ and  ${\rm Rep}_c(\sfG, \phi)''$ are equivalent tensor categories.
The consequences of Tannakian duality applied to C$^*$-algebras within  
$\cat_\sfG(\phi)={\rm Rep}_c(\sfG, \phi)$ will also be explored.

\begin{appendix}

\section{Equivalence to a strict category} \label{secappA}

MacLane showed that every monoidal category is equivalent to a strict category 
in which the structure maps such as $\Phi$ are the obvious identity maps.
The general construction is described and applied to the category 
$\cat_{\sfG}(\phi)$ in \cite{AM1, AM2}.
We shall give a simpler alternative construction which works in this case, and 
it is one of those cases where the structure is more transparent in the more 
general category $\cat^\sfH(\Phi)$ of modules for a Hopf algebra $\sfH$, although we 
shall apply it to $\sfH = C^*(\wh{{\sfG}}) \sim C_0({\sfG})$. (It is to some extent motivated 
by our observation that $\alg$-modules in the category $\cat_{\sfG}(\phi)$ are 
automatically modules for the crossed product $\alg\rtimes \wh{{\sfG}}$.) 

The algebra $\sfH$ is an $\sfH$-$\sfH$-bimodule under the left and right multiplication 
actions, and so, in particular it is an object in $\cat^\sfH(\Phi)$, though it is 
not an algebra in the category, since its multiplication is associative in the 
strict sense.
We shall exploit this dual role of $\sfH$ to construct the functor to a strict 
monoidal category.

Consider the functor $F$ which takes each object $\alg$ of $\cat^\sfH(\Phi)$ to 
the algebraic tensor product $F(\alg)= \alg\otimes \sfH$, and each morphism $T\in 
\hom(\alg,\bdd)$ to $F(T)\in \hom(F(\alg),F(\bdd))$ which sends $a\otimes h$ to 
$T(a)\otimes h$. 
(When $\sfH = C_0({\sfG})$ we can use $F(\alg) = C({\sfG},\alg)$ as a more convenient 
alternative.) 
Since $\sfH$ is a bimodule $F(\alg)$ is an $\sfH$-$\sfH$-bimodule, with 
the right multiplication action by $\sfH$, and the left comultiplication action, 
given in Sweedler notation by
$$h\cdot (a\otimes k) = \Delta h(a\otimes k) = h_{(1)}[a] \otimes h_{(2)}k.$$ 
With these actions we may take the new tensor product operation to be 
$\alg\otimes_F \bdd = \alg\otimes_\sfH \bdd$ (or $\alg\otimes_{C_0({\sfG})}\bdd$),  for 
which the new identity object is $F(\CC) = \sfH$ (or $C_0({\sfG})$) itself.

\begin{theorem} $\qquad$
\begin{itemize}
\item[(i)] For $\alg\in\cat^\sfH(\Phi)$ set $F(\alg) = \alg\otimes \sfH$ and  let 
the tensor product $F(\alg)\otimes_\sfH F(\bdd)$ be the quotient of 
$(\alg\otimes \sfH)\otimes(\bdd\otimes \sfH)$ by the equivalence relation that 
$(a\cdot h)\otimes b\sim a\otimes(h\cdot b)$, for all $a\in F(\alg)$, $b\in F(\bdd)$ and 
$h\in \sfH$. 
Then 
$$F(\alg)\otimes_\sfH F(\bdd) \cong F(\alg\otimes \bdd).$$ 

\item[(ii)] For $\alg \in \cat_{\sfG}(\phi)$ set $F(\alg) = C({\sfG},\alg)$ and let the tensor 
product $[F(\alg)\otimes_{C_0({\sfG})} F(\bdd)]$ be the quotient of 
$C({\sfG},\alg)\otimes C({\sfG},\bdd)$ by the equivalence relation that 
$(a\cdot h)\otimes b\sim a\otimes(h\cdot b)$, for all $a\in F(\alg)$, $b\in F(\bdd)$ and 
$h\in C_0({\sfG})$. 
Then 
$$F(\alg)\otimes_{C_0({\sfG})} F(\bdd) \cong F(\alg\otimes \bdd).$$ 
In each case set $F(T) = T\otimes{\rm id}$ for a morphism $T$.
Then $F$ defines a functor between tensor categories.
(The associator in each case is just $\Phi$ tensored with the identity.) 
\end{itemize}
\end{theorem}

\begin{proof}
We have
$$F(\alg)\otimes_\sfH F(\bdd) = (\alg\otimes \sfH)\otimes_\sfH (\bdd\otimes \sfH)
\cong (\alg\otimes \bdd)\otimes \sfH = F(\alg\otimes \bdd),$$
giving the result and showing consistency with the previous tensor product.
(In Sweedler notation we have the isomorphism is given explicitly by
$(a\otimes h)\otimes_H(b\otimes k)\mapsto  
(a\otimes \Delta(h)(b\otimes k)) =(a\otimes (h_{(1)}b)\otimes (h_{(2)}k)$.)
Most of the rest is easily checked.
In particular, we find that
\begin{align*}
((a\otimes h)\otimes_\sfH(b\otimes k))\otimes_\sfH(c\otimes l)
&= ((a\otimes h_{(1)}[b])\otimes (h_{(2)}k_{(1)}[c]))\otimes h_{(3)}k_{(2)}l\cr
&= \Phi[(a\otimes (h_{(1)}[b]\otimes h_{(2)}k_{(1)}[c]))]\otimes h_{(3)}k_{(2)}l\cr
&=(\Phi\otimes \id)[(a\otimes h)\otimes_\sfH[(b\otimes k)\otimes_\sfH(c\otimes l)]]\,.
\end{align*}
\end{proof}

For future reference we also note the connection with the crossed product.
 
\begin{proposition}
If $\alg$ is an algebra with multiplication $\star$ then $F(\alg)$ can be given 
the multiplication defined by the maps 
$$F(\alg)\otimes_\sfH F(\alg) = F(\alg\otimes \alg) \cong 
(\alg\otimes \alg)\otimes \sfH \longrightarrow (\alg\otimes \sfH) = F(\alg),$$ 
where the arrow denotes the map $\star\otimes 1$. 
When $\sfH = C_0({\sfG})$, $F(\alg)$ equipped with this 
multiplication is just the crossed product $\alg\rtimes {\sfG}$. 
If $\cM$ is a module for $\alg$ then $F(\cM)$ is a module for $F(\cA)$, and 
$({\rm id}\otimes \epsilon)F(\cM)$ can be identified with $\cM$ regarded as a 
$\alg\rtimes \wh{{\sfG}}$-module.
\end{proposition}

\begin{proof}
We have
$$(a\otimes h)\otimes(b\otimes k) = (a\otimes h_{(1)}[b])\otimes h_{(2)}k
\mapsto (a\star h_{(1)}[b])\otimes h_{(2)}k,$$
which is the crossed product multiplication in $\alg\rtimes {\sfG}$, when $\sfH = C_0({\sfG})$.
Replacing $b$ by an element of $\cM$  and applying $\epsilon$ (which is a counit 
and a multiplicative homomorphism) we have
\begin{align*}
({\rm id}\otimes \epsilon)[(a\otimes h)\star(m\otimes k)] 
&= (a\otimes h_{(1)}[m])\epsilon(h_{(2)}k)\\
&= (a\star h_{(1)}[m])\otimes\epsilon(h_{(2)})\epsilon(k)\\
&= (a\star h[m])\otimes\epsilon(k),
\end{align*}
as required.
\end{proof}

The advantage of expressing the crossed product action in terms of $F$ is that 
the functor respects direct sums and maps the free $\alg$-module $\alg^n$ to 
the  free $\alg\rtimes\wh{{\sfG}}$-module $(\alg\rtimes\wh{{\sfG}})^n$. 
If $\cM$ is a finite rank projective module defined by $e:\alg^n \to \cM$, 
then $F(\cM)$ is a finite rank projective module defined by 
$(\alg\rtimes\wh{{\sfG}})^n \to F(\cM)$, and, taking the image under 
${\rm id}\otimes \epsilon$, we see that $\cM$ is a finite rank projective 
$\alg\rtimes \wh{{\sfG}}$-module. 

We now introduce a new tensor product $F(\alg)\circ F(\bdd)$, such that for 
$a\in F(\alg)$, $b\in F(\bdd)$, and $c\in F(\cat)$, we have
$$(a\circ b)\otimes_\sfH c = a\otimes_\sfH (b\otimes_\sfH c).$$
This can be done explicitly using $\Phi$, 
since 
$$(a\circ b)\otimes_\sfH c = \Phi^{-1}((a\otimes_\sfH b)\otimes_\sfH c),$$
and $\Phi^{-1}$ is given by the left action of an element 
$\phi^{-1}\in \sfH\otimes \sfH\otimes \sfH$ on $\alg\otimes \bdd\otimes \cat$.
For convenience we shall write $\phi^{-1}$ in a Sweedler type notation as 
$\phi^{-1} 
= \phi^\prime\otimes\phi^{\prime\prime}\otimes\phi^{\prime\prime\prime}$, 
but this is neither intended to imply that $\phi$ is decomposable nor that this 
is in the range of 
$(\Delta\otimes 1)\Delta$. 
We then have $$(a\circ b)\otimes_\sfH c = 
(\phi^\prime\cdot a\otimes_\sfH \phi^{\prime\prime} \cdot b)\otimes_\sfH 
\phi^{\prime\prime\prime} c = (\phi^\prime\cdot a\otimes_\sfH 
\phi^{\prime\prime}\cdot b)\cdot \phi^{\prime\prime\prime}\otimes_\sfH c\,,$$ 
or, formally, 
$$a\circ b = (\phi^\prime\cdot a\otimes_\sfH 
\phi^{\prime\prime}\cdot b)\cdot \phi^{\prime\prime\prime} \in F(\alg)\otimes_\sfH F(\bdd).$$ 
In fact, this expression makes perfectly good sense (in the multiplier algebra 
if not in the original algebra), and can be used as a definition of $a\circ b$. 
(This would not have been true in the original setting with $a\in \alg$ and 
$b\in \bdd$, and the original tensor product, since there was only a left 
action and the argument would have led to 
$a\circ b  = (\phi^\prime\cdot a\otimes_\sfH 
\phi^{\prime\prime}\cdot b)\otimes\phi^{\prime\prime\prime} \in 
(\alg\otimes \bdd)\otimes \sfH$, not in $\alg\otimes \bdd$. 
Whether it makes sense in the original algebra or in some slightly extended 
algebra depends on the detail of the situation. 
The case when $\sfH=C_0({\sfG})$, $F(\alg) = C({\sfG},\alg)$ gives a large enough 
algebra to work whilst the algebraic tensor product would generally be too 
small.) 

Another useful way of thinking about the product is that 
$F(\alg\otimes \bdd) \sim [(\alg\otimes \bdd)\otimes \sfH]$, 
and $\Phi^{-1}$ maps this to $\alg\otimes (\bdd\otimes \sfH) 
\sim \alg\otimes F(\bdd)
\sim F(\alg)\otimes_\sfH F(\bdd)$.

As a consequence of the definition, for $d\in F({\cD})$
\begin{align*}
(a\circ b)\circ c)\otimes d 
&= (a\circ b)\otimes (c\otimes d)\cr
&= a\otimes(b\otimes(c\otimes d))\cr
&= a\otimes((b\circ c)\otimes d)\cr
&= (a\circ(b\circ c)) \otimes d,
\end{align*}
and taking ${\cD}=\sfH$ we see that we now have strict associativity.
We note that if $\sfH$ has a unit $1$ the original algebra $\alg$ can be 
identified with the subalgebra $\alg\otimes 1 \subseteq \alg\otimes \sfH$, 
(or in the case of $\sfH = C_0({\sfG})$ with the constant functions in the algebra 
$C({\sfG},\alg)$). 
These are not closed under the new tensor product $\circ$.

The new tensor product now carries over to products on algebra, so that 
$\star$ is replaced by a new product $a*b = \psi_z[a\star b]$, where $\phi_z(x,y) = \phi(x,y,z)$ 
acts as an element of $C({\sfG})\otimes C({\sfG})$,
The associativity of this may now be checked by either of the above calculations, and 
similarly for actions of algebras on modules.
In summary, the action of an algebra $\alg$ on a $C_0({\sfG})$-module can always be 
replaced by an action of the crossed product $\alg\rtimes \wh{{\sfG}}$, which a 
Fourier transform identifies with $C({\sfG},\alg)$, which means that for modules one 
always works with the objects whose multiplication can be made associative. 
This result clarifies  the paradox of the relevance of nonassociative 
algebras when algebras of operators are always associative, for we see that in 
representation the action of nonassociative algebras can always be replaced 
by an associative action if one so wishes. The situation is rather similar to 
that of projective representations of groups, where, for any multiplier 
$\sigma$ on a group ${\sfG}$, a projective $\sigma$-representation of ${\sfG}$ can be 
obtained from an ordinary representation of its central extension ${\sfG}^\sigma$. 
Nonetheless, there are situations in which  ${\sfG}$ and $\sigma$ appear naturally 
or linking a number of different situations, so that although ${\sfG}^\sigma$ is 
technically useful, it does not really capture the essence of the situation.
The study of the canonical commutation relations as a projective 
representation of a vector group provides a good example.
The central extension has representations allowing all possible values of 
Planck's constant, and so misses an important feature of the physical 
situation.

In addition to the left and right actions of $\sfH$, when $\sfH$ is a Hopf algebra 
there is also a right coaction of $\sfH$ on $F(\alg) = \alg\otimes \sfH$ given by 
$$\id\otimes \Delta: \alg\otimes \sfH \mapsto (\alg \otimes \sfH)\otimes \sfH,$$
and similarly for $C_0({\sfG})$.
There is no rebracketing problem since $\sfH$ is the algebra whose modules give 
the objects rather than an object itself, and the fact that this is a 
coaction follows from the coassociativity of $\Delta$.
To check compatibility with the tensor product structure, we note that
\begin{align*}
(a\otimes h)\otimes_\sfH(b\otimes k)
&= (a\otimes h_{(1)})\otimes_\sfH(b\otimes h_{(2)}k)\cr
&\mapsto [(a\otimes h_{(1)})\otimes_\sfH(b\otimes h_{(2)}k_{(1)})]\otimes h_{(3)}k_{(2)}\cr
&= [(a\otimes h_{(1)})\otimes_\sfH(b\otimes k_{(1)})]\otimes h_{(2)}k_{(2)},
\end{align*}
showing direct compatibility with the coaction on the individual factors.
(With a little modification of the functor $F$ it is possible to work with 
quasi-Hopf algebras too.)

If we instead consider the dual action of $\theta$ in the dual Hopf 
algebra $\sfH^*$ given in Sweedler notation by 
\begin{equation*}
\theta[a\otimes h] = \theta(h_{(2)})a\otimes h_{(1)} \,,
\end{equation*}
where $\theta(h_{(2)})$ denotes 
the pairing of $\sfH^*$ with $\sfH$, then this shows that one has the correct action 
on tensor products.
We can therefore regard $F$ as a functor from $\cat^\sfH(\Phi)$ to 
$\cat^{\sfH^*}(1)$ (or from $\cat_{\sfG}(\phi)$ to $\cat_{\wh{{\sfG}}}(1)$).

\subsection{Example: Strictification of the twisted compact operators}

The easiest way to see how the strictifcation works is to consider an example, 
such as the algebra $\cpt_\phi(L^2({\sfG}))$, which maps to 
$F(\cpt_\phi(L^2({\sfG}))) = C({\sfG},\cpt_\phi(L^2({\sfG})))$.
We first note that the product of $a\otimes h$ and $b\otimes k$ in 
$\alg\otimes C_0({\sfG})$ is given by
$$((a\otimes h)\star (b\otimes k))(x) = (a\star h_{(1)}[b])\otimes (h_{(2)}k)(x)
= (a\star h_{(1)}[b])\otimes h_{(2)}(x)k(x).$$
Now, defining $h_x(u) = h(ux) = (\Delta h)(u,x) = h_{(1)}(u)h_{(2)}(x)$,
we see that the product can be rewritten as
$$((a\otimes h)\star (b\otimes k))(x) = (a\star h_x[b])k(x).$$
Writing $k(x)(v,w) = k(v,w;x)$ for the $\cpt_\phi(L^2({\sfG}))$-valued function on 
${\sfG}$, the action of a function $h_x$ is just multiplication by 
$h_x(vw^{-1}) = h(vw^{-1}x)$
The product of two such kernel-valued functions is therefore 
\begin{equation*}
(k_1\star k_2)(u,w;x) 
= \int k_1(u,v;vw^{-1}x) k_2(v,w;x) \phi(uv^{-1},vw^{-1},w)\,dv\,.
\end{equation*}
This can be rewritten as
\begin{equation*}
(k_1\star k_2)(u,w;wx) 
= \int k_1(u,v;vx) k_2(v,w;wx)\phi(uv^{-1},vw^{-1},w)\,dv\,,
\end{equation*}
so that $k_j\mapsto k_j^\prime(u,w;x) = k_j(u,w;wx)$ gives an isomorphism with 
$C_0({\sfG})\otimes \cpt(L^2({\sfG}))$ equipped with componentwise multiplication.

According to the general prescription we can turn this into an associative 
product by applying $\Phi^{-1}$, that is, multiplying by 
$\phi(uv^{-1},vw^{-1},wx)^{-1}$, to get
\begin{align*}
&(k_1*k_2)(u,w;wx) \\
&= \int k_1(u,v;vx) k_2(v,w;wx)\phi(uv^{-1},vw^{-1},wx)^{-1}
\phi(uv^{-1},vw^{-1},w)\,dv \,,
\end{align*}
The cocycle identity tells us that 
\begin{equation*}
\phi(uw^{-1},w,x)\phi(uv^{-1},vw^{-1},wx) =
\phi(uv^{-1},vw^{-1},w)\phi(uv^{-1},v,x)\phi(vw^{-1},w,x)\,,
\end{equation*}
so that we can rewrite the product as
\begin{multline*}
(k_1* k_2)(u,w;wx)\phi(uw^{-1},w,x)^{-1} \\ 
= \int k_1(u,v;vx) k_2(v,w;wx) \phi(uv^{-1},v,x)^{-1} \phi(vw^{-1},w,x)^{-1}\,dv \,.
\end{multline*}
Setting $k^F(u,w;x) = k(u,w,wx)\phi(uw^{-1},w,x)^{-1}$ gives
\begin{equation*}
(k_1*k_2)^F(u,w,x)
= \int k_1^F(u,v;x) k_2^F(v,w;x)\,dv\,.
\end{equation*}
showing that the new product is isomorphic to the usual product on 
$C({\sfG},\cpt(L^2({\sfG})))$ \newline $\cong C_0({\sfG})\otimes \cpt(L^2({\sfG}))$,  which is certainly 
associative.

For future reference we note that the same argument applies to 
algebra-valued compact operators $\cpt_\phi(L^2({\sfG}))\otimes \alg$ when $\alg$ is 
associative, the only change being that the product of $k_j$ in the 
composition formula must be interpreted as a product in the algebra rather 
than of scalars.

This example shows that quite different nonassociative algebras can have the 
same associative version, since for any cocycle $\phi$, $\cpt_\phi(L^2({\sfG}))$ has 
$C({\sfG},\cpt(L^2({\sfG})))$ as its associative version (including the case of trivial 
$\phi$, when the algebra is associative from the start).

This is important in clarifying the duality theorem \cite[9.2]{BHM3}, where it was 
shown that the double dual of $\bdd = \uind{}_\sfN^{\sfG}(\alg)$ is 
$\cpt_\phi(L^2({\sfG}))\otimes \bdd$, since the associative version 
$\cpt(L^2({\sfG}))\otimes \bdd$ is stably equivalent to $\bdd$, as usual.
This certainly means that all the representations (that is, the modules) 
correspond naturally with those of $\bdd$.
We summarise this in a theorem.

\begin{theorem}
The dual $\wh{\bdd} = \bdd\rtimes_{\beta,v}{\sfG}$ of the algebra 
$\bdd = \uind_N^{\sfG}(\alg)$ is a nonassociative algebra with a natural action of 
$\wh{{\sfG}}$, and the double dual $\wh{\bdd}\rtimes \wh{{\sfG}}$ can be given the 
associative product $\cpt(L^2({\sfG}))\otimes \uind_N^{\sfG}(\alg)$.
\end{theorem}

We conclude by noting that the octonions $\OO$ can be constructed 
from $\RR$ using ${\sfG} = \ZZ_2\times \ZZ_2\times \ZZ_2$ with the 
cocycle
\begin{equation}
\phi({\bf a},{\bf b},{\bf c}) = (-1)^{[{\bf a},{\bf b},{\bf c}]}
= (-1)^{{\bf a}\cdot ({\bf b}\times {\bf c})}\,,
\end{equation}
where ${\sfG}$ is identified with $\{0,1\}^3\subset \RR^3$.
Thus the same procedure can be used to give an associative version of 
$C(\ZZ_2\times \ZZ_2\times \ZZ_2,\OO)$, as a 64-dimensional 
algebra over $\RR$.

\section{Nonassociative bounded operators, tempered \\ distributions \&
a concrete approach to nonassociative $C^*$-algebras} \label{secappB}

We begin with an illustrative example. Let $\sfG=\RR^n$, and consider the space of all 
bounded operators  $\cb(L^2(\sfG))$ on the Hilbert space $L^2(\sfG)$. We begin by 
showing that $T \in \cb(L^2(\sfG))$ determines a unique tempered distribution $k_T$ 
on $\sfG^2$.  That is, there is a canonical embedding,
$\cb(L^2(\sfG)) \hookrightarrow \cS'(\sfG^2)$.
This embedding will be frequently used, for instance to give the algebra
$\cb(L^2(\sfG))$ a nonassociative product, which has the advantage of 
being rather explicit. Later, we will also determine other closely related results. 
Recall
that the Sobolev spaces $H^s(\sfG), \, s\in \RR$, 
are defined as follows: 
the Fourier transform on Schwartz functions 
on $\sfG$ is a topological isomorphism,
$\,\,\,\what{}\, \colon \cS(\sfG) \to   \cS(\sfG),$
where we identify $\sfG$ with its Pontryagin dual group. 
It extends uniquely to an isometry on square integrable 
functions on $\sfG$,
$\,\,\,\what{} \ \colon L^2(\sfG) \to   L^2(\sfG)$.
Moreover, by duality, the Fourier transform extends to be a topological isomorphism
on tempered distributions on $\sfG$,
$\,\,\,\what{} \ \colon \cS'(\sfG) \to   \cS'(\sfG)$.
Then for $s\in \RR$, define $H^s(\sfG)$ to be the Hilbert space of all 
tempered distributions $Q$ such that $(1+|\xi|^2)^{s/2} \what Q(\xi)$ is in $L^2(\sfG)$, 
with inner product $\langle Q_1, Q_2\rangle_s = 
\langle (1+|\xi|^2)^{s/2}\what Q_1(\xi), (1+|\xi|^2)^{s/2}\what Q_2(\xi)\rangle_0$, where 
$\langle\ ,\ \rangle_0$ denotes the inner product on $L^2(\sfG)$.

The following are some basic properties of Sobolev spaces,
which are established in any basic reference on distribution theory. 
For $s<t$, $H^t(\sfG) \subset H^s(\sfG)$ and moreover the inclusion map 
$H^t(\sfG) \hookrightarrow H^s(\sfG)$ is continuous. Also one has 
$\cS(\sfG) = \bigcap_{s\in \RR} H^s(\sfG)$, $\,\,\cS'(\sfG) = \bigcup_{s\in \RR} H^s(\sfG)$
and the inclusions $\iota_s \colon \cS(\sfG)   \hookrightarrow H^s(\sfG)$
and $\kappa_s \colon H^s(\sfG) \hookrightarrow \cS'(\sfG) $ are continuous for any $s\in \RR$.
The renowned Schwartz kernel theorem says that a continuous linear operator
$T \colon \cS(\sfG) \to \cS'(\sfG)$ determines a unique tempered distribution $k_T$ 
on $\sfG^2$, and conversely. 

\begin{lemma} There is a canonical embedding,
\begin{equation}
\cb(L^2(\sfG)) \hookrightarrow \cS'(\sfG^2)\,,
\end{equation}
whose image is contained in the subspace of composable tempered distributions.
\end{lemma}

\begin{proof} Suppose  that $T \in \cb(L^2(\sfG))$. Then in the notation above, 
the composition
\begin{equation}
\kappa_0 \circ T \circ \iota_0 \colon \cS(\sfG) \to \cS'(\sfG) \,,
\end{equation}
is a continuous linear operator. By the Schwartz kernel theorem, it determines a 
unique tempered distribution $k_T \in  \cS'(\sfG^2)$. Suppose now that $S 
\in \cb(L^2(\sfG))$. Then $ST \in  \cb(L^2(\sfG))$ and 
\begin{equation}
k_{ST}(x,y) = \int_{z\in \sfG} k_S(x, z) k_T(z, y)\ dz \,,
\end{equation}
where  $\displaystyle \int_{z\in \sfG}\, dz$ denotes the distributional pairing.  
\end{proof}

We can now define a new product on $\cb(L^2(\sfG))$ making it into a nonassociative 
C$^*$-algebra.

\begin{definition}
Let $\phi \in C(\sfG\times \sfG\times \sfG)$ be an antisymmetric tricharacter on $\sfG$.
For $S, T \in \cb(L^2(\sfG))$, define the tempered distribution $k_{S\star T} \in  \cS'(\sfG^2)$
by the formula
\begin{equation}
k_{S\star T}(x,y) = \int_{z\in \sfG} k_S(x, z) k_T(z, y) \phi(x,y,z)\ dz\,.
\end{equation}
Then for all $\xi, \psi\in L^2(\sfG)$,  the linear operator
$S\star T$ given by the prescription
\begin{equation}
\langle \xi, S\star T \psi \rangle_0 =  \int_{x, y\in \sfG} k_{S\star T}(x,y) \bar \xi(x) \psi(y)\ dxdy\,,
\end{equation}
defines a bounded linear operator  in $\cb(L^2(\sfG))$, which follows from the earlier observation
that $S\star T $ is an adjointable operator.

This extends the definition in \cite{BHM3} of twisted compact operators $\cK_\phi(L^2(\sfG))$. 
Then by \S4, $\,\,\star$ defines a nonassociative product on $\cb(L^2(\sfG))$ 
which agrees with the nonassociative product on the twisted compact operators, which will be 
justified 
in what follows. We denote by $\cb_\phi(L^2(\sfG))$ the space $\cb(L^2(\sfG))$ endowed with the 
nonassociative product $\star$, and call it the algebra of twisted bounded operators.
\end{definition}

There is an involution $k_{S^*}(x,y) = \conj{k_S(y,x)}$ for all $S\in \cb_\phi(L^2(\sfG))$, 
and the norm on $\cb_\phi(L^2(\sfG))$ is the usual operator norm.
The following are obvious from the definition:
$\, \forall \; \lambda \in \CC,\; \forall \; S_1, S_2 \in \cb_\phi(L^2(\sfG))$, 
\begin{equation}\label{B*identities1}
\begin{array}{rcl}
(S_1+S_2)^* & = &  S_1^* + S_2^*\,,\\
(\lambda S_1)^* & = &  \bar\lambda S_1^*\,,\\
S_1^{**} & = &  S_1\,.\\
\end{array}
\end{equation}
The following lemma can be proved as in Section 5 in \cite{BHM3}.

\begin{lemma}  $\forall \; S_1, S_2 \in \cb_\phi(L^2(\sfG))$, 
\begin{equation}\label{B*identities2}
(S_1\star S_2)^* = S_2^* \star S_1^*\,.
\end{equation}
\end{lemma}

What appears to be missing for the deformed bounded operators  $\cb_\phi(L^2(\sfG))$
is the so called C$^*$-identity,
\begin{equation}\label{C*identities}
||S_1^* \star S_1|| = ||S_1^* S_1|| = ||S_1||^2.
\end{equation}
However, we will continue to call $\cb_\phi(L^2(\sfG))$ a nonassociative C$^*$-algebra
and this prompts the following definition of a general class of nonassociative C$^*$-algebras.

\begin{definition}
A nonassociative $C^*$-subalgebra $\cA$ of $\cb_\phi(L^2(\sfG))$, 
is defined to be a $\sfG$-invariant, $\star$-subalgebra of $\cb_\phi(L^2(\sfG))$ that 
is closed under taking adjoints and also closed in the operator norm topology.
\end{definition}

In particular, such an $\cA$ satisfies the identities in equations \eqref{B*identities1} and \eqref{B*identities2}.
Examples include the algebra of twisted bounded operators $\cb_\phi(L^2(\sfG))$ and 
the algebra of twisted compact operators $\cK_\phi(L^2(\sfG))$.
The following two propositions can be proved as in Section 5 of \cite{BHM3}.

\begin{proposition}
The group $\sfG$ acts on the twisted algebra of bounded operators $\cb_\phi(L^2(\sfG))$ 
by natural $*$-automorphisms
\begin{equation}
\theta_x[k](z,w) = \phi(x,z,w)k(zx,wx),
\end{equation}
and $\theta_x\theta_y = \ad(\sigma(x,y))\theta_{xy}$,
where $\ad(\sigma(x,y))[k](z,w) = \phi(x,y,z)k(z,w)\phi(x,y,w)^{-1}$ comes from 
the multiplier $\sigma(x,y)(v) = \phi(x,y,v)$.
\end{proposition}

\begin{proposition}
$\cb_\phi(L^2(\sfG))$ is a continuous 
deformation of $\cb(L^2(\sfG))$.
\end{proposition}

\section{Nonassociative crossed products and nonassociative tori} \label{secappC}

\subsection{Nonassociative tori -- revisited} \label{secappCa}

Here will present a slightly different, more geometric, approach to the definition of the 
nonassociative torus as defined in \cite{BHM3}, which in fact generalizes 
the construction there, and also realizes it as a nonassociative 
deformation of the algebra continuous functions on the torus.
We begin with a general construction, and later specialize to the case when
$M$ is the torus.

{\bf Basic Setup:}
{\em Let $M$ be a compact manifold with fundamental group $\Gamma$, 
and $\wM$ its universal cover. Assume 
for simplicity that $\wM$ is contractible, that is $M=B\Gamma$ is the classifying space
of $\Gamma$.  In that case we have an isomorphism $H^n(M,\ZZ) \cong H^n(\Gamma,\ZZ)$,
due to Eilenberg and MacLane \cite{EMacLa, EMacLb}.}

A large class of examples 
of manifolds $M$ that satisfy the hypotheses of the Basic Setup are locally symmetric 
spaces $M = \Gamma\backslash \sfG/\sfK$, where $\sfG$ is a Lie group, $\sfK$ a maximal compact 
subgroup of $\sfG$, $\Gamma$ a discrete, torsion-free cocompact subgroup of $\sfG$, 
since in this case $\wM = \sfG/\sfK$ is a contractible manifold. 
This includes tori and hyperbolic manifolds in particular. 

The isomorphism $H^n(M,\RR) \to H^n(\Gamma,\RR)$ can be explicitly constructed 
by making use of the double complex $( C^{p,q}(\wM,\Gamma) ; \delta,d)$, with
\begin{equation}
C^{p,q}(\wM,\Gamma) = C^p(\Gamma, \Omega^q(\wM)) =
\{ f : \Gamma^{\otimes p} \to \Omega^q(\wM) \} \,,
\end{equation}
where we think of $q$-forms on $\wM$ as a left $\Gamma$-module through the 
action $\gamma \cdot \omega = (\gamma^*)^{-1} \omega$.
The differential $d : C^{p,q}(\wM,\Gamma) \to C^{p,q+1}(\wM,\Gamma)$ is the de Rham 
differential on $\Omega(\wM)$, hence this complex is acyclic since $\wM$ is contractible.
The differential $\delta : C^{p,q}(\wM,\Gamma) \to C^{p+1,q}(\wM,\Gamma)$ is given by
\begin{equation}
(\delta f)_{\gamma_1,\ldots,\gamma_{p+1}} = \gamma_1 \cdot f_{\gamma_2,\ldots,\gamma_{p+1}}
- f_{\gamma_1\gamma_2,\ldots,\gamma_{p+1}} + (-1)^p f_{\gamma_1,\ldots,\gamma_p\gamma_{p+1}}
+ (-1)^{p+1} f_{\gamma_1,\ldots,\gamma_p} \,.
\end{equation}
Its cohomology $H^p(\Gamma, \Omega^q(\wM))$ is known as the group cohomology with
coefficients in the $\Gamma$-module $\Omega(\wM)$.  
The map $H(M,\RR) \to H(\Gamma,\RR)$ is now obtained by `zigzagging' through this
double complex, much in the same way as for the \v Cech-de Rham complex.  We will 
illustrate the procedure in the case of our interest, i.e.\ how to explicitly associate 
to a closed degree 3 differential form $H$ on $M$ a ${\mathsf{U}(1)}$-valued 
3-cocycle $\sigma$ on the discrete group $\Gamma$. 

First we let $\wH$ denote the  lift of $H$ to $\wM$.  
Now $\wH = dB$, where $B$ is a 2-form on $\wM$, i.e.\
$B \in C^{0,2}(\wM,\Gamma)$. Since $\wH$ is $\Gamma$-invariant, it follows that 
for all $\gamma \in \Gamma$, we have $0 = \gamma \cdot \wH - \wH = d (\gamma \cdot B - B) =
d\delta B$, 
so that $\gamma\cdot B - B$ is a closed 2-form on $\wM$.  By hypothesis, 
it follows that  
\begin{equation}\label{eqn:H1}
(\delta B)_\gamma = \gamma\cdot B - B = dA_\gamma\,,
\end{equation}
for some 1-form $A_\gamma$ on $\wM$, i.e.\ $A\in C^{1,1}(\wM,\Gamma)$. 
Then by \eqref{eqn:H1}, it follows that the following identity holds for all
$\beta , \gamma \in \Gamma$:
\begin{equation*} 
d( \beta\cdot A_\gamma - A_{\beta\gamma} + A_\beta ) = d \delta A = \delta dA = 
\delta^2 B = 0 \,.
\end{equation*}
By hypothesis, it follows that  
\begin{equation*}
\beta\cdot A_\gamma - A_{\beta\gamma} + A_\beta= ( \delta A)_{\beta,\gamma} = df_{\beta, \gamma} 
\end{equation*}
for some smooth function $f_{\beta, \gamma} $ on $\wM$, that is,\
$f\in C^{2,0} (\wM,\Gamma)$.  Continuing, one computes that,
\begin{equation*}
d(\alpha\cdot f_{\beta, \gamma} - f_{\alpha\beta, \gamma} + f_{\alpha, \beta\gamma} - 
f_{\alpha, \beta}) =  d \delta f = \delta^2 A = 0\,.
\end{equation*}
Therefore 
\begin{equation*}
\alpha\cdot f_{\beta, \gamma} - f_{\alpha\beta, \gamma}  + f_{\alpha, \beta\gamma} - f_{\alpha, \beta} 
=(\delta f )_{\alpha,\beta,\gamma} = c(\alpha, \beta,\gamma) \,,
\end{equation*}
where $c(\alpha, \beta, \gamma)$ is a constant.  That is,\ 
$c: \Gamma \times \Gamma \times \Gamma \to \RR$
is a 3-cocycle on $\Gamma$, and we can set for all $t\in \RR$,
\begin{equation}
\sigma_t(\alpha, \beta, \gamma) = \exp\left(it c(\alpha, \beta, \gamma)\right)\,.
\end{equation}	
Then $\sigma_t(\alpha, \beta, \gamma)$ is a ${\sf U}(1)$-valued 3-cochain on $\gamma$, 
satisfying the pentagonal identity,
\begin{equation}
\sigma_t(\alpha, \beta, \gamma) \sigma_t(\alpha, \beta\gamma, \delta) 
\sigma_t(\beta, \gamma, \delta)
= \sigma_t(\alpha\beta, \gamma, \delta) \sigma_t(\alpha, \beta, \gamma\delta)\,,
\end{equation}
for all $\alpha, \beta, \gamma, \delta \in \Gamma$.
That is, $\sigma_t(\alpha, \beta, \gamma) $ is actually a ${\sf U}(1)$-valued 
3-cocycle on $\Gamma$.

It is convenient to normalize the function $f_{\alpha, \beta}$ such that 
$f_{\alpha, \beta}(x_0) = 0$ for all $\alpha, \beta \in \Gamma$ and for some $x_0 \in \wM$.
Then the formula for the  ${\sf U}(1)$-valued 3-cocycle on $\Gamma$ simplifies to,
\begin{equation*}
\sigma_t(\alpha, \beta, \gamma) = \exp\left(it(\alpha^{*-1} f_{\beta, \gamma}(x_0))\right)\,.
\end{equation*}	

Consider the unitary operator $u_t(\beta, \gamma)$ acting on $L^2(\Gamma)$ given by
\begin{equation*}
(u_t(\beta, \gamma) \psi)(\alpha) = \exp(it \alpha^{*-1} f_{\beta, \gamma}(x_0) )\psi (\alpha) = 
\sigma_t(\alpha, \beta, \gamma) \psi(\alpha).
\end{equation*}
We easily see that 
\begin{equation*}
\sigma_t(\alpha, \beta, \gamma)  u_t(\alpha, \beta) u_t(\alpha\beta, \gamma)  = 
\xi_\alpha[u_t(\beta, \gamma)]  u_t(\alpha, \beta\gamma) 
\end{equation*}
where $\xi_\alpha = {\rm ad}(\rho(\alpha))$ and $(\rho(\alpha)\psi)(g) = \psi(g\alpha)$
is the right regular representation. 
One can define, analogous to what was done in \cite{BHM3}, a twisted convolution product and 
adjoint on $C(\Gamma, \cK)$, $\cK = \cK(L^2(\Gamma))$,  by 
\begin{equation}
(f*_t g)(x) = \sum_{y\in \Gamma} f(y) \xi_y[g(y^{-1}x)] u_t(y, y^{-1}x)  \,,
\end{equation}
and 
\begin{equation}
f^*(x) = u_t(x, x^{-1})^{-1} \xi_x[f(x^{-1})]^* \,.
\end{equation}
The operator norm completion is the nonassociative twisted crossed product C$^*$-algebra 
\begin{equation*}
C^*(\cK, \Gamma, \sigma_t) = \cK(L^2(\Gamma)) \rtimes_{\xi, u_t} \Gamma.
\end{equation*} 
where as before, $\sigma_t(\alpha, \beta, \gamma)$ is a ${\sf U}(1)$-valued 3-cocycle on $\Gamma$
as above.
This construction extends easily to the case when $\cK$ is replaced by a general 
$\Gamma$-C$^*$-algebra $\cA$, giving rise to a nonassociative C$^*$-algebra denoted by 
$C^*(\cA, \Gamma, \sigma_t) = \cA \rtimes_{\xi, u_t} \Gamma.$ 

In the special case when $M=\TT^n$ is a torus, we get the nonassociative torus $A_{\sigma_t}(n)$. 
Now $A_{\sigma_0}(n)$ is just the ordinary crossed product $\cK(L^2(\ZZ^n)) \rtimes \ZZ^n$,
where $\ZZ^n$ acts on $\cK(L^2(\ZZ^n)) $ via the the adjoint of the left regular representation. 
By the Stabilization theorem, and using the Fourier transform,
\begin{equation*}
\cK(L^2(\ZZ^n)) \rtimes \ZZ^n \cong C^*(\ZZ^n) \otimes \cK(L^2(\ZZ^n)) \cong C(\TT^n) \otimes 
 \cK(L^2(\ZZ^n)).
\end{equation*}
This then indicates why $A_{\sigma_t}(n)$ is a nonassociative deformation of the ordinary torus 
$\TT^n$ for $t\ne 0$.  

\begin{example}
As an explicit example, consider the 3-torus $M=\TT^3$.  We have $\wM=\RR^3$ and $\Gamma=
\ZZ^3$.   Let us take, for $H\in H^3(M,\RR)$, $k$ times the volume form on $M$
(i.e.\ $k$ times the image in de Rham cohomology of the generator 
of $H^3(M,\ZZ)\cong \ZZ$).  Its lift to $\wM$ is explicitly given by
\begin{equation*}
\wH = k \, dx_1 \wedge dx_2 \wedge dx_3 \,,
\end{equation*}
where $(x_1,x_2,x_3)$ are standard coordinates on $\wM = \RR^3$.   
Let us denote elements of $\Gamma=\ZZ^3$ by $\mathbf n = (n_1,n_2,n_3)$.  
Going through the procedure above, we see that a representative of this 
3-form in group cohomology is given by $c(\mathbf l,\mathbf m,\mathbf n) =k\,  l_1m_2n_3$.
However, by making different choices for $B, A_{\mathbf n}$, etc., 
specifically
\begin{align*}
B & = \txtfrac13 k\, 
  ( x_1\, dx_2 \wedge dx_3 + \text{cycl} ) \,,\\
A_{\mathbf n} & =  \txtfrac16 k\, \left( n_1 ( x_2 dx_3 - x_3 dx_2) + \text{cycl} \right )  \,, \\
f_{\mathbf m, \mathbf n} & = \txtfrac16 k\, \left( m_1 ( n_2 x_3 - n_3 x_2)+ \text{cycl} \right )  \,,
\end{align*}
we can also construct
a completely antisymmetric representative, namely
\begin{equation}
c(\mathbf l,\mathbf m,\mathbf n) = \txtfrac16 k\,  \mathbf l \cdot (\mathbf m \times \mathbf n) \,.
\end{equation}
It is this representative which gives rise to an antisymmetric tricharacter $\sigma_t$ on
$\Gamma$.  Note that the image of an integer cohomology class in $H^3(M,\RR)$
is not necessarily integer-value, but in this example rather ends up in $\frac{1}{6}\ZZ $.  
This example explains and corrects a discrepancy between our earlier paper \cite{BHM3}
and the physical interpretation of our nonassociative 3-torus in the context of
open string theory by Ellwood and Hashimoto (cf.\ Eqn. (5.14) in \cite{EH}).
\end{example}

\subsection{Factors of automorphy and continuous trace C$^*$-algebras} \label{secappCb}

We next study principal $\PU$-bundles $P$ and associated bundles of compact operators $\cK_P$
and their sections over manifolds $M$ that satisfy the assumptions of the Basic Setup of
Appendix \ref{secappCa}. 
Let $\wP$ denote the lift of $P$ to $\wM$.
Since $H^3(\wM)=0$, it follows that $\wP$ is trivializable, i.e. $\wP\cong \wM \times \PU$. 
Having fixed this isomorphism, we can define a continuous map 
$j : \Gamma \times \wM \to \PU = {\rm Aut}(\cK)$ by the following commutative diagram,
\begin{equation}
\xymatrix{
                 \cK = (\cK_\wP)_x  \ar[ddr]_{p} \ar[rr]^{j(\gamma, x)} && 
                 \cK = (\cK_\wP)_{\gamma\cdot x}   \ar[ddl]^p                   \\  && \\ 
             &        (\cK_P)_{p(x)} &      }
\end{equation}
Then 
\begin{equation}\label{eqn:auto}
j(\gamma_1\gamma_2, x) ^{-1}j(\gamma_1, \gamma_2 x) j(\gamma_2, x) = 1.
\end{equation}
is a factor of automorphy for the bundle $\cK_P$. Conversely, given a continuous map 
$j : \Gamma \times \wM \to \PU = {\rm Aut}(\cK)$ satisfying \eqref{eqn:auto}, we can define
a bundle of compact operators, 
\begin{equation}
\cK_{j} = (\wM \times \cK)/\Gamma \to M
\end{equation}
where $\gamma \cdot (x, \xi) = (\gamma\cdot x, j(\gamma, x) \xi)$ for $\gamma \in \Gamma$ and 
$(x, \xi) \in \wM \times \cK$. 

Given any two algebra bundles of compact operators, $\cK_j, \cK_{j'}$ over $M$
with factors of automorphy $j, j'$ respectively, and an isomorphism $\phi : \cK_j \longrightarrow
\cK_{j'}$, we get an induced isomorphism 
\begin{equation}
\widetilde \phi : \wM \times \cK =  \widetilde\cK_j 
\longrightarrow \cK_{j'} =  \wM \times \cK
\end{equation}
given by $\widetilde\phi(x, \xi) = (x, u(x)\xi)$, 
where $u: \wM \to \PU$ is continuous. Since $\widetilde\phi$ commutes with the action of $\Gamma$,
$\gamma\cdot  \widetilde\phi(x, \xi) = \widetilde\phi(\gamma(x, \xi))$, we deduce that 
$
(\gamma\cdot x, j'(\gamma, x) u(x)\xi) = (\gamma\cdot x, u(\gamma\cdot x)j(\gamma, x)\xi).
$
Therefore 
\begin{equation}\label{eqn:autoequiv} 
 j'(\gamma, x) = u(\gamma\cdot x)j(\gamma, x) u(x)^{-1}
 \end{equation}
for all $x\in \wM$ and 
$\gamma \in \Gamma$. Conversely, two factors of automorphy $j, j'$ give rise to 
isomorphic algebra bundles $\cK_j, \cK_{j'}$ of compact operators if they
are related by \eqref{eqn:autoequiv} for some continuous function $u: \wM \to \PU$.

Therefore continuous sections of $\cK_P$ can be viewed as continuous maps $f\in C(\wM, \cK)$ satisfying 
the property, 
\begin{equation}\label{C14}
f(\gamma\cdot  x) = j(\gamma, x) f(x),\qquad \forall \gamma\in \Gamma, \, x\in \wM.
\end{equation}
For example, $f(x):=\sum_{\gamma \in \Gamma} j(\gamma, x)^{-1} F(\gamma \cdot x)$
converges uniformly on compact subsets of $\wM$ whenever $F: \wM \to \cK$ is a compactly supported continuous function, and satisfies \eqref{C14}, therefore defining a continuous section of $\cK_P$.

We would like to write the Dixmier-Douady invariant of the algebra bundle of compact operators 
$\cK_P$ over $M$ in terms of the factors of automorphy.
There is no obstruction to lifting the factor of automorphy $j : \Gamma \times \wM \to \PU = 
{\rm Aut}(\cK)$
to $\wh j: \Gamma \times \wM \to \U$, because of our assumptions on $\wM$. However 
the cocycle condition \eqref{eqn:auto} has to be modified, 
\begin{equation}
\wh j(\gamma_1\gamma_2, x) ^{-1} \wh j(\gamma_1, \gamma_2 x) \wh j(\gamma_2, x) = 
\tau(\gamma_1, \gamma_2, x)\,,
\end{equation}
where $\tau : \Gamma \times \Gamma \times \wM \to \U(1)$.
There is no obstruction to lifting $\tau: \Gamma \times\Gamma \times\wM \to \U(1)$ 
to $\wh\tau :  \Gamma \times\Gamma \times\wM \to \RR$, however the cocycle condition
satified by $\tau$ has to be modified to $\delta\wh\tau(\gamma_1, \gamma_2, \gamma_3) = 
\eta(\gamma_1, \gamma_2, \gamma_3)$, where $\eta:  \Gamma \times\Gamma \times \Gamma
\to \ZZ$ is a $\ZZ$-valued 3-cocycle on $\Gamma$. One can show  $DD(P) = \delta([\tau']) = [\eta]$.
Thus, given a principal $\PU$ bundle $P$ on $M$, we have derived a cohomology 
class $ [\eta] \in H^3(\Gamma, \ZZ) \cong H^3(M, \ZZ)$ which is by standard arguments 
independent of the choices made. The relation with the previous discussion is that $[\eta] = [c]$. 

To see that the converse is true, notice that $\tau$ can be viewed as a 
continuous map $\tau': \Gamma\times\Gamma \to C(\wM, \U(1))$,
which is easily verified to be a $C(\wM, \U(1))$-valued 2-cocycle on $\Gamma$. 
Recall from standard group cohomology theory that equivalence 
classes of extensions of a group $\Gamma$ by an abelian 
group $C(\wM, \U(1))$ on which $\Gamma$ acts is in bijective correspondence
correspondence with the group cohomology with coefficients, 
$H^2(\Gamma, C(\wM, \U(1)))$. 
We will first show that possible extensions $\wGamma$ of $\Gamma$ by $C(\wM, \U(1))$
are in bijective correspondence with elements of $H^3(M, \ZZ)$ called the Dixmier-Douady 
invariant, and we will also compute $DD(P) \in H^3(M, \ZZ)$ in our case. 
Now there is an exact sequence of abelian groups,
\begin{equation}
0 \to \ZZ \to C(\wM, \RR) \to C(\wM, \U(1)) \to 0 \,.
\end{equation}
This leads to a change of coefficients long exact sequence,
\begin{equation}
\cdots \to H^2(\Gamma,  C(\wM, \RR)) \to H^2(\Gamma, C(\wM, \U(1))) \stackrel{\delta}{\to} 
H^3(\Gamma, \ZZ) \to 
H^3(\Gamma,  C(\wM, \RR))  \to \cdot
\end{equation}
Since $\Gamma$ acts freely on $\wM$ and $C(\wM, \RR)$ is an induced module, it follows that 
$ H^j(\Gamma,  C(\wM, \RR)) = 0$ for all $j>0$.  Therefore $H^j(\Gamma, C(\wM, \U(1)))
\cong H^{j+1}(\Gamma, \ZZ) = H^{j+1}(M, \ZZ)$ for all $j>0$, and in particular for $j=2$ as claimed.
In particular, since $[\tau'] \in H^2(\Gamma, C(\wM, \U(1)))$, we see that $DD(P) = \delta([\tau'])=[\eta]
\in H^3(\Gamma, \ZZ) = H^3(M, \ZZ)$ is the Dixmier-Douady invariant of $P$.
 
We next explain how the data $(H, B, A_\gamma, f_{\alpha, \beta})$ also 
determine a bundle gerbe in a natural way.
The bundle gerbe consists of the submersion $\wM \to M$. Then the fibered product $\wM^{[2]}$ is equivariantly 
isomorphic to $\Gamma \times \wM$. Under this identification, the two projection maps 
$\pi_i : \wM^{[2]} \to \wM, \, i=1, 2$ become the action $\mu : \Gamma \times \wM \to \wM$ of 
$\Gamma$ on $\wM$ and the projection $p: \wM \times \Gamma \to \wM$, respectively. 
Then $A_\gamma$ defines a 
connection on the trivial line bundle $\cL_\gamma \to \{\gamma\} \times  \wM $ whose curvature 
is $dA_\gamma$. The choice of curving is $B$, satisfying the equation 
$\mu^*B - p^* B = dA$, which on the sheet $\{\gamma\} \times  \wM $ reduces to 
$\gamma^* B - B = dA_\gamma$.
The 3-curvature $dB = \wH$ is the lift of the closed 3-form $H$ on $M$. What is surprising is that 
$H$ need not have integral periods!


\section{Motivation from T-duality in String Theory} \label{secappD}

For completeness we have summarized the original motivation for this 
work, namely T-duality in string theory, in this appendix.  
 We believe, however, that the results in
this paper are of interest independent of our original motivation.

T-duality, also known as target space duality, plays an important role in string theory
and has been the subject of intense study for many years. 
In its most basic form, T-duality relates a string theory compactified
on a circle of radius $R$, to a string theory compactified on the dual circle of radius $1/R$
by the interchange of the string momentum and winding numbers.
T-duality can be generalized to locally defined circles (principal circle bundles, circle fibrations),
higher rank torus bundles or fibrations, and in the presence of a background H-flux which is
represented by a closed, integral \v Cech 3-cocycle $H$ on the spacetime manifold 
$Y$.  It is closely related to mirror symmetry 
through the SYZ-mechanism.

An amazing feature of T-duality is that it can relate topologically distinct spacetime manifolds by the
interchange of topological characteristic classes with components of the H-flux.  
Specifically, denoting by $(Y,[H])$ the pair of an (isomorphism class of) principal 
circle bundle $\pi: Y\to X$, characterized by the first Chern class $[F]\in H^2(X,\ZZ)$ 
of its associated line bundle,  and an H-flux $[H] \in H^3(Y,\ZZ)$, the T-dual again 
turns out to be a pair $(\widehat Y,[\widehat H])$, where the principal circle bundles

$$\xymatrix{
\TT  \ar[r] &  Y \ar[d]_{\pi}  \\
& X}, \qquad \qquad \qquad \qquad
\xymatrix{
\TT  \ar[r] &  \widehat Y \ar[d]_{\widehat\pi}  \\
& X}
$$
are related by
$[\widehat F] = \pi_*[H] ,\  [F] = \widehat \pi_* [\widehat H]$,
such that on the correspondence space
$$
\xymatrix @=4pc @ur { Y \ar[d]_{\pi} & 
Y\times_X  \widehat Y \ar[d]_{\widehat p} \ar[l]^{p} \\ X & \widehat Y\ar[l]^{\widehat \pi}}
$$
we have $p^* [H] - \widehat p^*[\widehat H] = 0$ \cite{BEM03a,BEM03b}. 

In earlier papers we have argued that the twisted K-theory $K^\bullet (Y,[H])$ 
(see, e.g., \cite{BCMMS01})
classifies charges of D-branes on $Y$ in the background of H-flux $[H]$ \cite{BM00},
and indeed, as a consistency check, one can prove that T-duality gives an isomorphism 
of twisted K-theory (and the closely related twisted cohomology $H^\bullet(Y,[H])$ by 
means of the twisted Chern character $ch_H$) \cite{BEM03a}
$$\xymatrix{
K^\bullet(Y, [H])  \ar[r]^{T_!} \ar[d]^{ch_H} &  
      K^{\bullet+1}(\widehat Y, [ \widehat H]) \ar[d]^{ch_{\widehat H} } \\         
H^\bullet  (Y, [H])    \ar[r]^{T_*} &   H^{\bullet+1} (\widehat Y, [\widehat H])     
}
$$
The above considerations were generalized to principal torus bundles in
\cite{BHM04a,BHM04b}.

Since the projective unitary group of an infinite dimensional Hilbert space
$\PUH$ is a model for $K(\ZZ,2)$, we can `geometrize' the H-flux in terms of an
(isomorphism class of) principal $\PUH$-bundle $P$ over $Y$.  We can reformulate the 
discussion of T-duality above in terms of noncommutative geometry as follows.
The space of continuous sections vanishing at infinity, $\alg = C_0(Y,\mathcal E)$, 
of the associated algebra bundle of compact operator $\cK$ on the Hilbert space 
$\mathcal E = P \times_{\PUH} \mathcal K$. $\mathcal A$
is a stable, continuous-trace, C$^*$-algebra
with spectrum $Y$, and has the property that it is locally Morita equivalent to continuous 
functions on $Y$. Thus the $H$-flux has the effect of making spacetime noncommutative.
The K-theory of $\alg$ is just the twisted K-theory $K^\bullet (Y,[H])$.
The $\TT$-action on $Y$ lifts essentially uniquely to an $\RR$-action
on $\alg$.  In this context T-duality is the operation of taking the
crossed product $\alg \rtimes \RR$, which turns out to be another continuous 
trace algebra associated to $(\widehat Y, [\widehat H])$ as above. A fundamental property of T-duality is that 
when applied twice, yields the original algebra $\alg$, and the reason that it works in this
case is due to Takai duality.   The isomorphism of the D-brane charges 
in twisted K-theory is, in this context, due to the Connes-Thom isomorphism.
These methods have been generalized to principal torus bundles
by Mathai and Rosenberg \cite{MR04a, MR04b, MR05}, however novel features arise.
First of all the $\TT^n$-action on the principal
torus bundle $Y$ need not always lift to an $\RR^n$-action on
$\alg$.  Even if
it does, this lift need not be unique.  Secondly, the crossed product $\alg \rtimes \RR^n$
need not be a continuous-trace algebra, but rather, it could be a continuous field of 
noncommutative tori 
\cite{MR04a}, and necessary and sufficient conditions are given when these T-duals occur.  More generally, 
as argued in \cite{BHM04a}, when the $\TT^n$-action on the principal
torus bundle $Y$ does not lift to an $\RR^n$-action on
$\alg$, one has to leave the category of $C^*$-algebras
in order to be able to define a ``twisted" lift. The associator in this case is the restriction 
of the H-flux $H$ to the torus fibre of $Y$, and the ``twisted" crossed product is defined
to be the T-dual. The fibres of the T-dual are noncommutative,  nonassociative tori. 
That this is a proper definition of T-duality is due to our results which show that 
the analogs of Takai duality and the Connes-Thom isomorphism hold in this context.
Thus an appropriate context to describe nonassociative algebras that arise as T-duals 
of spacetimes with background flux, such as nonassociative tori, is that of C$^*$-algebras 
in tensor categories, which is the subject of this paper.

\end{appendix}

\end{document}